\documentclass[final,leqno]{siamltex}
\usepackage{mathbbol}
\usepackage{mathrsfs}
\usepackage{amsfonts}
\usepackage{}

\usepackage{stmaryrd}
\usepackage{amsmath}
\usepackage{amssymb}
\usepackage{multirow}
\usepackage{cite}
\usepackage{caption}
\usepackage{tabularx}
\usepackage{stmaryrd}
\usepackage{booktabs}
\usepackage{cite}
\usepackage{amssymb}
\usepackage{verbatim}
\usepackage[all,cmtip]{xy}
\usepackage{color}
\usepackage{graphicx}
\usepackage{subfigure}

\input psfig.sty

\newcommand{\p}{\partial}
\newcommand{\og}{\omega}
\newcommand{\Og}{\Omega}
\newcommand{\ep}{\varepsilon}
\newcommand{\fl}[2]{\frac{#1}{#2}}
\newcommand{\dt}{\delta}
\newcommand{\nn}{\nonumber}

\newcommand{\tht}{\theta}

\newcommand{\wt}{\widetilde}
\newcommand{\Dt}{\Delta}
\renewcommand{\theequation}{\arabic{section}.\arabic{equation}}
\newcommand{\be}{\begin{equation}}
\newcommand{\ee}{\end{equation}}
\newcommand{\ba}{\begin{array}}
\newcommand{\ea}{\end{array}}
\newcommand{\bea}{\begin{eqnarray}}
\newcommand{\eea}{\end{eqnarray}}
\newcommand{\beas}{\begin{eqnarray*}}
\newcommand{\eeas}{\end{eqnarray*}}
\newtheorem{remark}{Remark}[section]
\newcommand{\bx}{{\bf x} }

\title{Uniform error bounds of a finite difference method for the Zakharov
system in the subsonic limit regime via an asymptotic
consistent formulation\thanks{This work was partially supported by the Ministry
of Education of Singapore grant
R-146-000-196-112 (W. Bao) and the Natural Science Foundation
of China Grant 91430103 (C. Su).}}

\author{Weizhu Bao\thanks{Department of Mathematics,
National University of Singapore, Singapore 119076 ({\tt matbaowz@nus.edu.sg},
URL: http://www.math.nus.edu.sg/\~{}bao/)} \and
Chunmei Su\thanks{Beijing Computational Science Research Center, Beijing 100193,
China; and The Fields Institute for Research in Mathematical Sciences,
222 College Street, University of Toronto, Toronto, Ontario M5T 3J1, Canada ({\tt sucm@csrc.ac.cn})}}

\date{}
\begin{document}

\maketitle

\begin{abstract}
We present a uniformly accurate finite difference method and
establish rigorously its uniform error bounds for
the Zakharov system (ZS) with a dimensionless parameter $0<\varepsilon\le 1$,
which is inversely proportional to the speed of sound.
In the subsonic limit regime, i.e., $0<\varepsilon\ll 1$, the solution propagates
highly oscillatory waves and/or rapid outgoing initial layers due to the perturbation of the
wave operator in ZS and/or the incompatibility of the initial data which
is characterized by two nonnegative parameters $\alpha$ and $\beta$. Specifically, the solution
propagates waves with $O(\varepsilon)$- and $O(1)$-wavelength in time and space, respectively,
and amplitude at $O(\varepsilon^{\min\{2,\alpha,1+\beta\}})$ and
$O(\varepsilon^\alpha)$ for well-prepared ($\alpha\ge1$) and ill-prepared ($0\le \alpha<1$)
initial data, respectively. This high oscillation of the solution in time brings
significant difficulties in designing numerical methods and
establishing their error bounds, especially in the subsonic limit regime.
A uniformly accurate finite difference method is proposed by reformulating ZS into an
asymptotic consistent formulation and adopting an integral approximation of
the oscillatory term. By adapting the energy method and using the limiting equation
via a nonlinear Schr\"{o}dinger equation with an oscillatory potential,
we rigorously establish two independent error bounds at
$O(h^2+\tau^2/\varepsilon)$ and $O(h^2+\tau^2+
\tau\varepsilon^{\alpha^*}+\varepsilon^{1+\alpha^*})$, respectively, with $h$ the mesh size,
$\tau$ the time step and $\alpha^*=\min\{1,\alpha\}$.
Thus we obtain error bounds at
$O(h^2+\tau^{4/3})$ and $O(h^2+\tau^{1+\frac{\alpha}{2+\alpha}})$
for well-prepared and ill-prepared initial data,
respectively, which are uniform in both space and time for
$0<\varepsilon\le 1$ and optimal at the second order in space. Other techniques in the analysis include
the cut-off technique for treating the nonlinearity and inverse estimates
to bound the numerical solution.
Numerical results are reported to demonstrate that our error bounds are sharp.
\end{abstract}


\pagestyle{myheadings}\thispagestyle{plain}

\section{Introduction}
\setcounter{equation}{0}
\label{section1}
Consider the dimensionless Zakharov system (ZS) for describing the propagation of
Langmuir waves in plasma \cite{Sulem,Ozawa}
\be\label{Zak}
\begin{split}
&i\p_tE^\ep(\bx,t)+\Dt E^\ep(\bx,t)-N^\ep(\bx,t) E^\ep(\bx,t)=0,  \quad\bx \in \mathbb{R}^d, \quad t>0,\\
&\ep^2 \p_{tt} N^\ep(\bx,t)-\Dt N^\ep(\bx,t)-\Dt |E^\ep(\bx,t)|^2=0,  \quad\bx \in \mathbb{R}^d, \quad t>0,\\
&E^\ep(\bx,0)=E_0(\bx),\quad N^\ep(\bx,0)=N_0^\ep(\bx),\quad \p_t N^\ep (\bx,0)=N_1^\ep(\bx), \quad\bx\in \mathbb{R}^d.
\end{split}
\ee
Here $t$ is time, $\bx$ is the spatial coordinates,
the complex function $E^\ep:=E^\ep(\bx,t)$ is the slowly varying
envelope of the highly oscillatory electric field, the real function
$N^\ep:=N^\ep(\bx,t)$ represents the deviation of the ion density from its equilibrium value,
$0<\ep\le 1$ is a dimensionless parameter which is inversely proportional to the acoustic speed,
and $E_0(\bx)$, $N_0^\ep(\bx)$ and  $N_1^\ep(\bx)$ are given functions satisfying $\int_{\mathbb{R}^d}
N_1^\ep(\bx)d\bx=0$.

There exist extensive analytical and numerical studies in the literatures for
the standard ZS, i.e. $\ep=1$ in \eqref{Zak}.  Along the analytical
part, for the derivation of ZS from the Euler-Poisson equations,
we refer to \cite{Ginibre, Sulem}; and
for the well-posedness, we refer to \cite{Bourgain, Colliander, Ginibre, Sulem}
and references therein. Based on these results, we know that the ZS \eqref{Zak}
conserves the {\sl wave energy}
\be\label{mass}
\mathcal{M}(t)=
\|E^\ep(\cdot,t)\|^2:=\int_{\mathbb{R}^d}|E^\ep(\bx,t)|^2d\bx\equiv
\int_{\mathbb{R}^d}|E_0(\bx)|^2d\bx=\mathcal{M}(0), \quad t\ge0,
\ee
and the {\sl Hamiltonian}
\begin{equation} \label{energy}
\mathcal{L}^\ep(t):=\int_{\mathbb{R}^d}\left[|\nabla E^\ep|^2+N^\ep|E^\ep|^2
+\fl{1}{2}\left(\ep^2|\nabla U^\ep|^2+|N^\ep|^2\right)\right]d\bx
\equiv \mathcal{L}^\ep(0),\quad t\ge0,
\end{equation}
where $U^\ep:=U^\ep(\bx,t)$ is defined as
\be\label{Udef}
-\Delta U^\ep(\bx,t)=\partial_t N^\ep(\bx,t), \quad \bx\in \mathbb{R}^d, \qquad
\lim\limits_{|\bx|\rightarrow \infty} U^\ep(\bx,t)=0, \qquad t\ge0.
\ee
Along the numerical part, different numerical methods
have been proposed and analyzed in the last two
decades.  Glassey \cite{Glassey} presented an
energy-preserving implicit finite difference scheme and established
an error bound at first order in both spatial and temporal discretizations.
Later, Chang and Jiang \cite{Chang} improved it to the optimal second order
convergence by considering an implicit or semi-explicit conservative
finite difference schemes \cite{Chang2}. Other approaches include the
exponential-wave-integrator spectral method \cite{Bao2003, Payne}, Jacobi-type method
\cite{Bhrawy}, Legendre-Galerkin method \cite{Ji},
discontinuous-Galerkin method \cite{Xia} and time-splitting spectral
method \cite{Bao2005, Jin}. The analytical and numerical results
for ZS have been extended to
the generalized Zakharov system \cite{Hadouaj, Hadouaj2}, the
vector Zakharov system \cite{Sulemb} and the vector Zakharov
system for multicomponents \cite{Hadouaj2}.

When $\ep\to0^+$, i.e., in the subsonic limit regime, formally
we get $E^\ep(\bx,t)\to E(\bx,t)$, $\rho^\ep:=\rho^\ep(\bx,t)=|E^\ep|^2\to |E|^2=\rho$ and $N^\ep(\bx,t)\to N(\bx,t)=-|E(\bx,t)|^2$, where $E:=E(\bx,t)$ satisfies
the cubic nonlinear Schr\"{o}dinger equation (NLSE)
\cite{Masmoudi2008, Ozawa,Schochet}
\be\label{NLS}
\begin{split}
&i\p_t E(\bx,t)+\Dt E(\bx,t)+|E(\bx,t)|^2E(\bx,t)=0,\quad t>0,\quad \bx \in \mathbb{R}^d,\\
& E(\bx,0)=E_0(\bx),\quad \bx \in \mathbb{R}^d.
\end{split}
\ee
The NLSE (\ref{NLS}) conserves the wave energy (\ref{mass}) with $E^\ep=E$
and the {\sl Hamiltonian}
\be\label{energy-nls}
\mathcal{L}(t):=\int_{\mathbb{R}^d}\left[|\nabla E(\bx,t)|^2-\fl{1}{2}|E(\bx,t)|^4\right]d\bx
\equiv\mathcal{L}(0), \quad t\ge0.
\ee
Convergence rates of the subsonic limit from the ZS (\ref{Zak})
to the NLSE (\ref{NLS}) and initial layers as
well as the propagation of oscillatory waves have been rigorously studied in
the literatures \cite{Masmoudi2008, Ozawa,Schochet}. Based on the results,
when $0<\ep\ll1$, the solution of the ZS (\ref{Zak}) propagates
highly oscillatory waves at wavelength $O(\ep)$ and $O(1)$ in time and space, respectively,
and/or rapid outgoing initial layers at speed $O(1/\ep)$ in space.
In addition, the initial data ($E_0, N_0^\ep, N_1^\ep$) in \eqref{Zak}
can be decomposed as
\bea\label{initZS}
\begin{split}
&\qquad \ N_0^\ep(\bx)=N(\bx,0)+\ep^\alpha \og_0(\bx),\quad
N_1^\ep(\bx)=\partial_tN(\bx,0)+\ep^\beta \og_1(\bx),\qquad \bx \in \mathbb{R}^d,\\
&\qquad \ N(\bx,0)=-|E_0(\bx)|^2, \quad \partial_tN(\bx,0)=-\partial_t\rho(\bx,0)=2\mathrm{Im}(\Dt E_0(\bx)\overline{E_0(\bx)}):=\phi_1(\bx),
\end{split}
\eea
where $\alpha, \beta\ge0$ are parameters describing the
incompatibility of the initial data of the ZS \eqref{Zak} with respect to
 that of the NLSE \eqref{NLS} in the subsonic limit regime, 
$\og_0(\bx)$ and $\og_1(\bx)$
are two given real functions independent of $\ep$ and satisfy $\int_{\mathbb{R}^d}\og_1(\bx)d\bx=0$, and
$\mathrm{Im}(f)$ and $\overline{f}$ denote the imaginary and complex conjugate parts of $f$,
respectively.
In fact, when $\alpha\ge2$ and $\beta\ge1$,
the leading order oscillation is due to the term $\ep^2\partial_{tt}N$ in ZS;
and when either $0\le \alpha<2$ or $0\le \beta<1$, the leading order oscillation is due to
the initial data.

\begin{figure}[t!]
\begin{minipage}[t]{0.5\linewidth}
\centering
\includegraphics[width=2.6in]{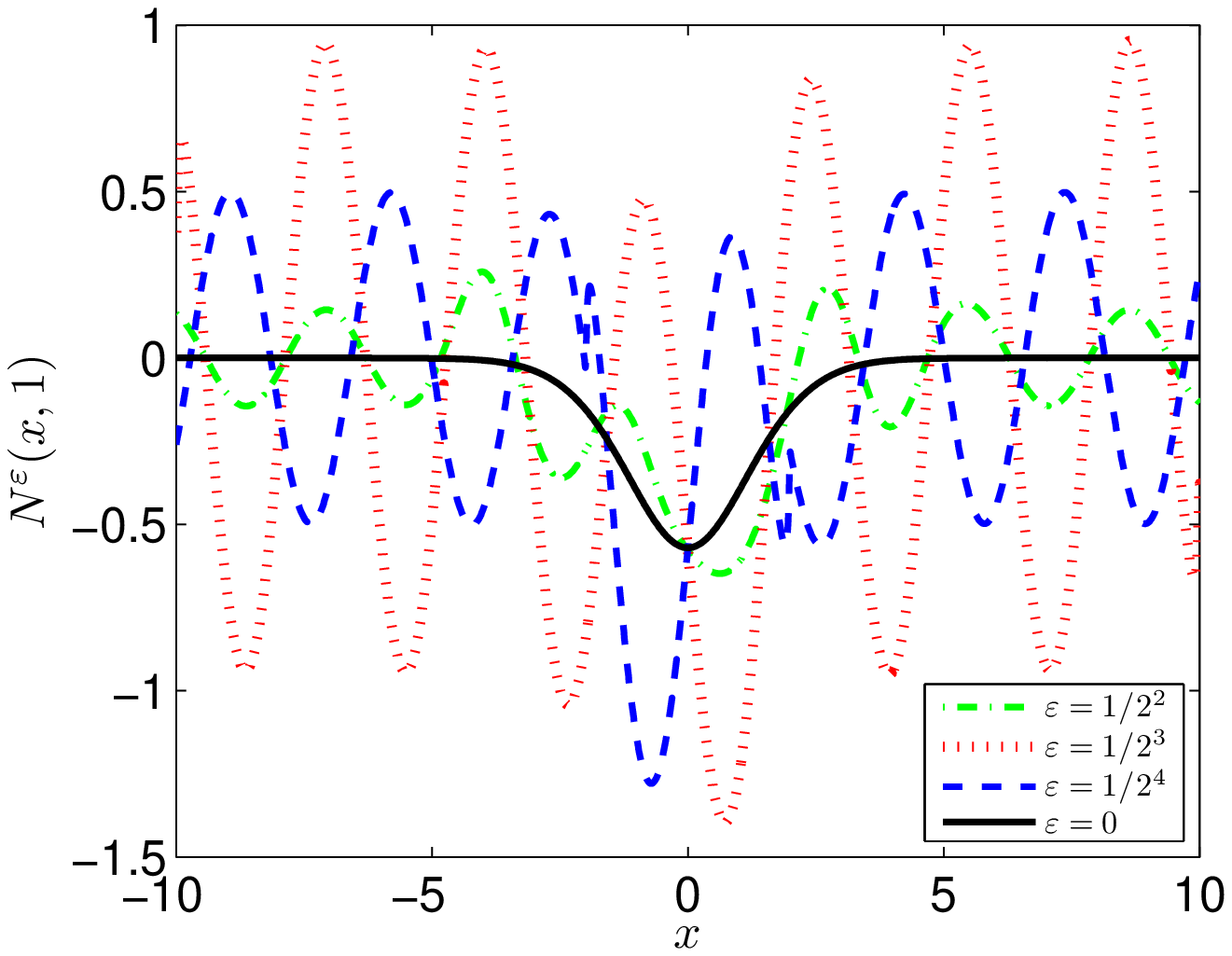}
\end{minipage}%
\begin{minipage}[t]{0.5\linewidth}
\centering
\includegraphics[width=2.6in]{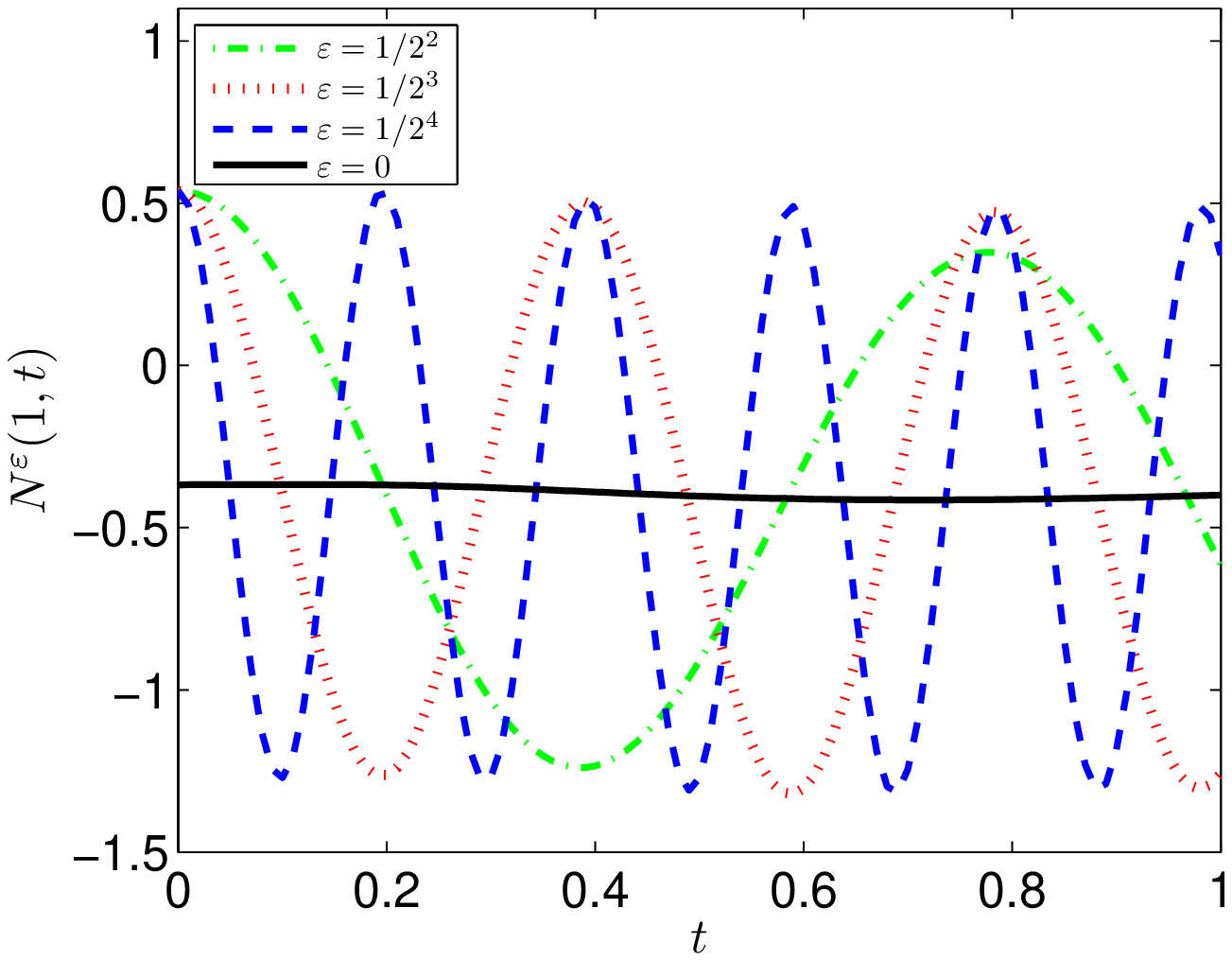}
\end{minipage}
\begin{minipage}[t]{0.5\linewidth}
\centering
\includegraphics[width=2.6in]{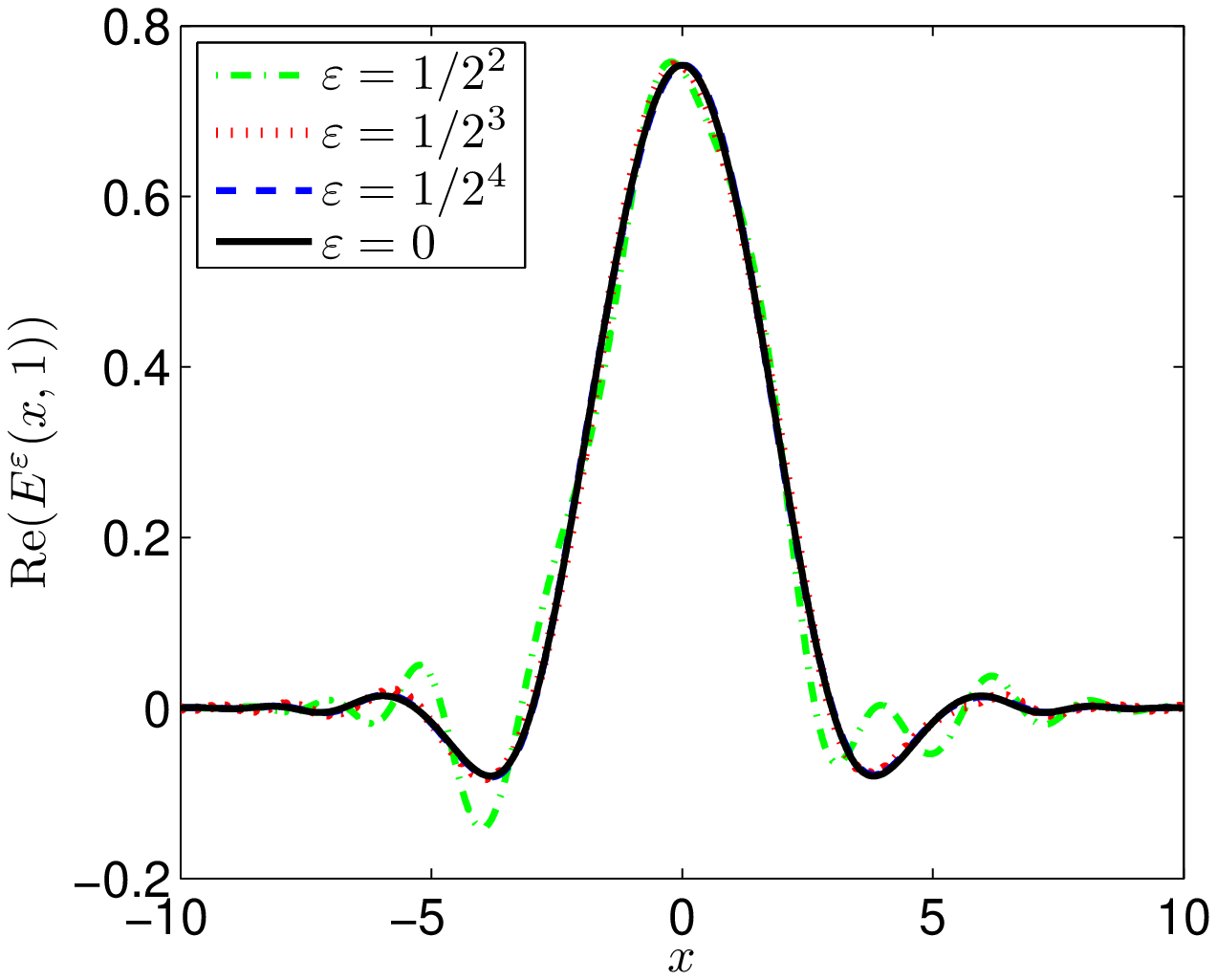}
\end{minipage}%
\begin{minipage}[t]{0.5\linewidth}
\centering
\includegraphics[width=2.6in]{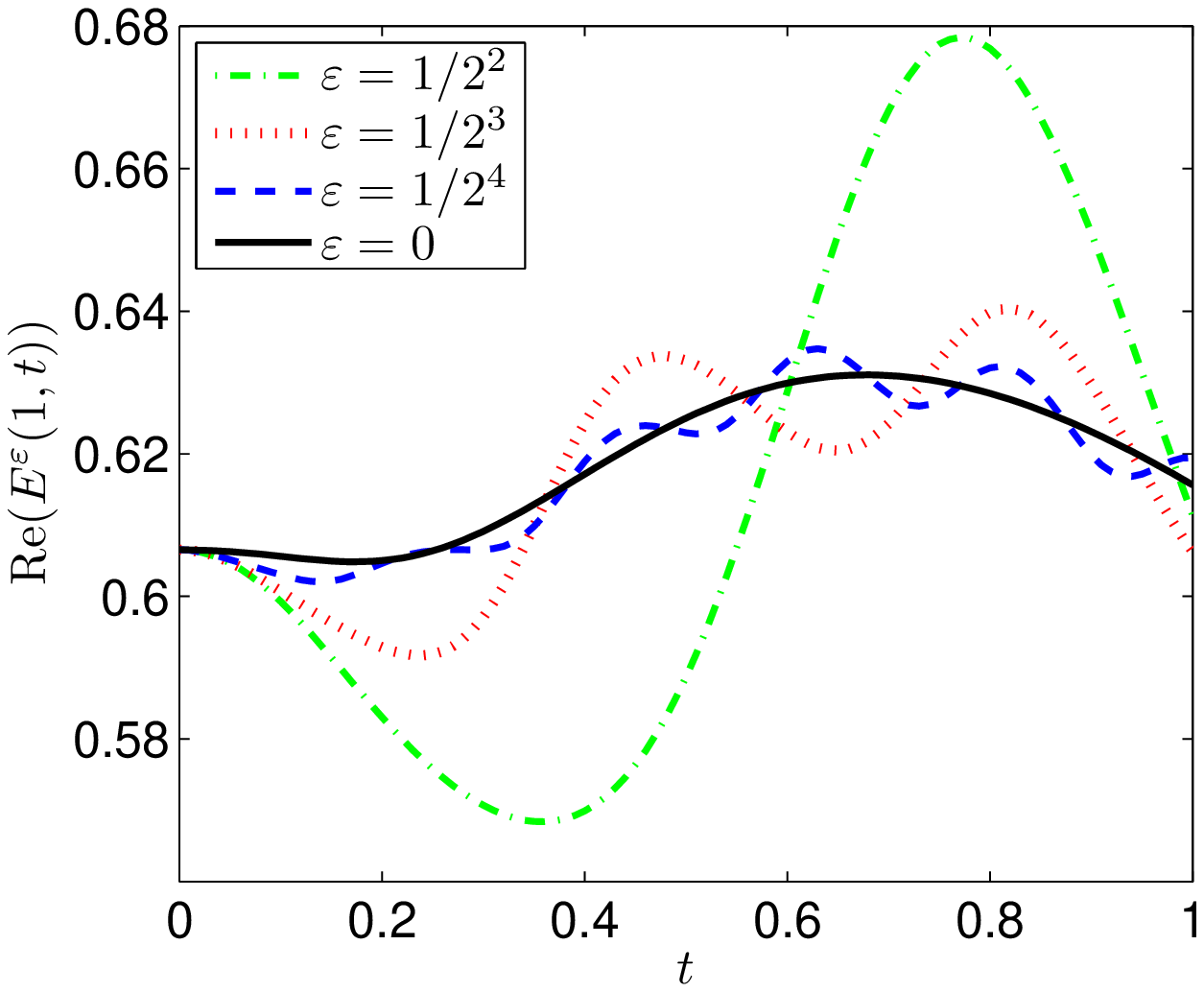}
\end{minipage}
\begin{minipage}[t]{0.5\linewidth}
\centering
\includegraphics[width=2.6in]{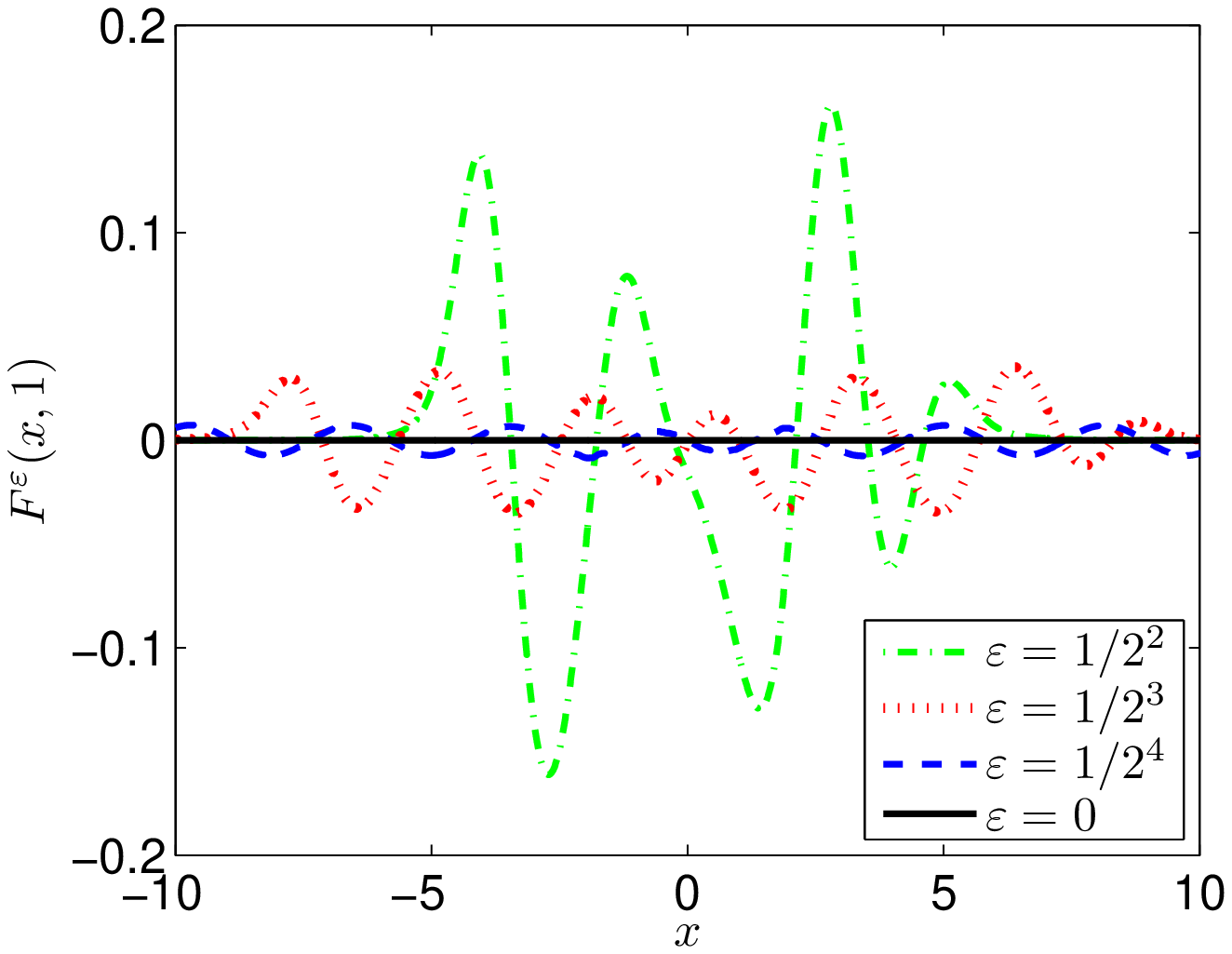}
\end{minipage}%
\begin{minipage}[t]{0.5\linewidth}
\centering
\includegraphics[width=2.6in]{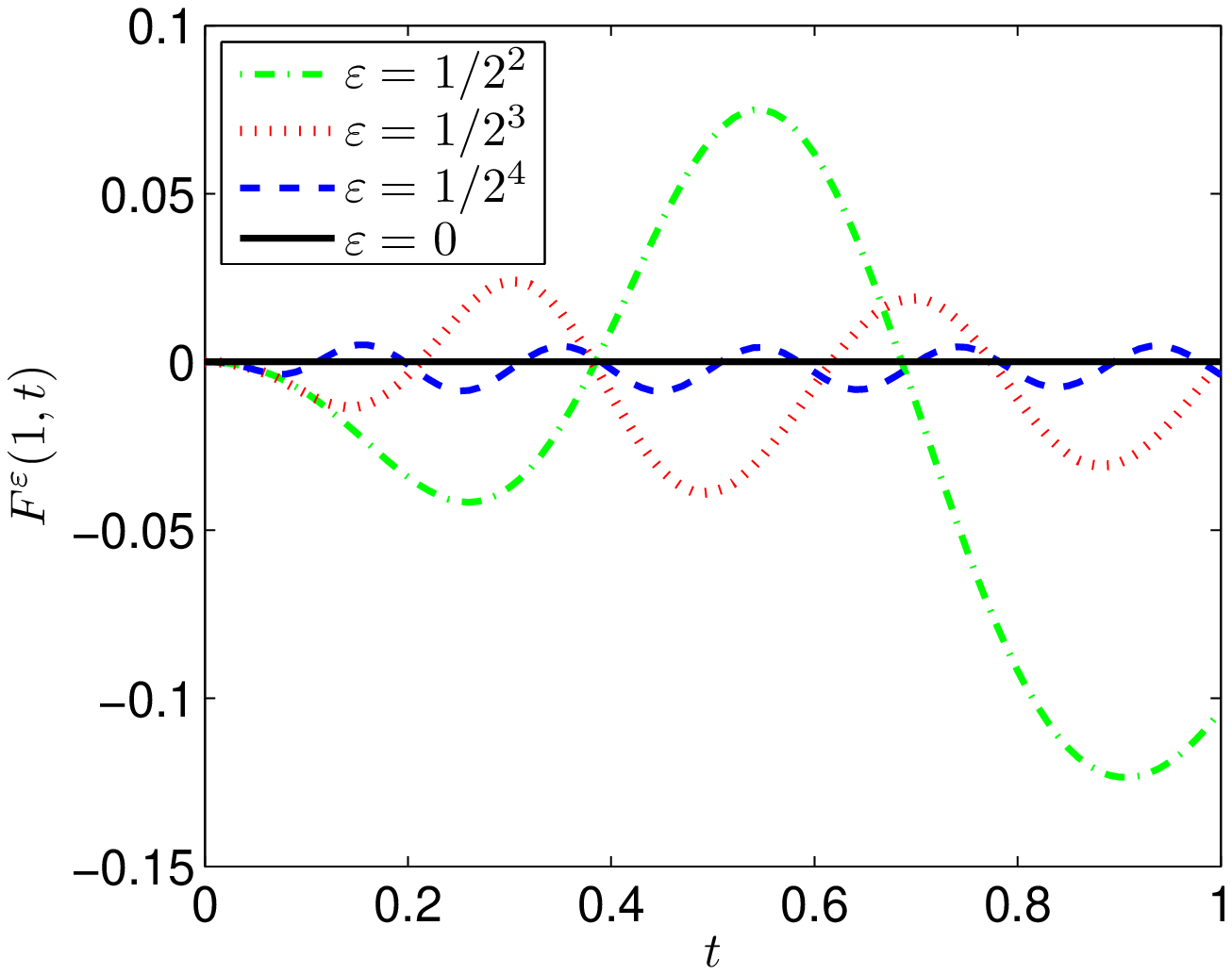}
\end{minipage}\caption{The solutions of the ZS \eqref{Zak} for
different $\ep>0$ and the NLSE ($\ep=0$) as well as
$F^\ep$ defined in \eqref{Nd} with $d=1$. Here $\mathrm{Re}(f)$
denotes the real part of $f$.}\label{fig1}
\end{figure}

To illustrate the oscillatory and/or rapid outgoing wave phenomena,
Fig. \ref{fig1} shows the solutions $N^\ep(x,1)$, $N^\ep(1,t)$,
Re$(E^\ep(x,1))$ and  Re$(E^\ep(1,t))$ of the
ZS \eqref{Zak} with $d=1$, $E_0(x)=e^{-x^2/2}$,
$\alpha=0$, $\beta=0$, $\og_0(x)=e^{-\fl{1}{18^2-x^2}}\sin (2x)\chi_{(-18,18)}$
with $\chi$  the characteristic function and
 $\og_1(x)\equiv0$ in (\ref{initZS}) for different $\ep$, which was
obtained numerically on a bounded computational interval $[-200,200]$ with
the homogenous Dirichlet boundary condition \cite{Bao2005}. For comparison,
here we also plot  $F^\ep(x,1)$ and $F^\ep(1,t)$  defined in \eqref{Nd}.

The highly oscillatory nature of the solution of the ZS (\ref{Zak}) in time
brings significant numerical burdens, especially in the subsonic limit regime.
Some numerical results for ZS with different $0<\ep\le 1$ have been
reported in the literatures \cite{Bao2005,Jin}. To the best of our knowledge,
there are few results concerning error estimates
of different numerical methods for ZS with respect to
the mesh size $h$, time step $\tau$ as well as the parameter $0<\ep\le1$
except that an error bound of the finite difference Legendre
pseduospectral method was derived for ZS in one dimension (1D)
when $\alpha\ge2$ and $\beta\ge1$ \cite{Ji}.
Very recently, for the conservative finite difference method,
Cai and Yuan \cite{Cai2015} established uniform error bounds at $O(h^2+\tau^{4/3})$
for $0<\ep\le1$ when $\alpha \ge2$ and $\beta \ge 1$, and at  $O(h^2+\tau^{\fl{2}{3}\min\{\alpha,1+\beta\}})$
when $1\le \alpha <2$ and/or $0\le \beta<1$.
However, when $0<\alpha<1$, their error bound $O(h^2/\ep^{1-\alpha}+\tau^{\fl{2}{3}\alpha})$ is not uniform
in space, and in particular, when $\alpha=0$, their error bound
$O(h^2/\ep+\tau^2/\ep^3)$ requests the meshing strategy (or $\ep$-scalability)
$h=O(\ep^{1/2})$ and $\tau=O(\ep^{3/2})$ which is not uniform in both space and time
when $0<\ep\ll1$. The reason is due to that
$N^\ep(\bx,t)$ does not converge to  $N(\bx,t)=-|E(\bx,t)|^2$ when $\alpha=0$
and $\ep\to0^+$ \cite{Masmoudi2008,Schochet,Sulemb} (cf. Fig 1.1 top row).

The aim of this work is to design a finite difference method
for ZS, which is uniformly accurate in space and time for $0<\ep\le1$,
and carry out rigorous error analysis for the finite difference method by
paying particular attention to how the error bounds depend on explicitly
$h$ and $\tau$ as well as the parameter $\ep$. The key ingredients in
designing the uniformly accurate finite difference method are based on
(i) reformulating ZS into an asymptotic consistent formulation and (ii) adapting
an integral approximation of the oscillatory term. In establishing error bounds,
we adapt the energy method, cut-off technique for treating the nonlinearity,
the inverse estimates to bound the numerical solution, and the limiting equation
via a nonlinear Schr\"{o}dinger equation with an oscillatory potential.
The error bounds of our new numerical method significantly improve
the results of the standard finite difference method for ZS in the subsonic limit regime
\cite{Cai2015}, especially for the ill-prepared initial data, i.e. $0\le \alpha< 1$.

The rest of the paper is organized as follows. In section 2, we introduce
an asymptotic consistent formulation of ZS, present a finite difference method
and state our main results. Section 3 is devoted to the details of the error analysis.
Numerical results are reported in section 4 to confirm our error bounds. Finally
some conclusions are drawn in section 5. Throughout the paper, we adopt the standard
Sobolev spaces and the corresponding norms and adopt $A\lesssim B$ to mean that
there exists a generic constant $C>0$ independent of $\ep$, $\tau$, $h$,
such that $|A|\le C\,B$.

\section{A finite difference method and its error bounds}
In this section, we will introduce
an asymptotic consistent formulation of ZS, present a uniformly accurate
finite difference method and state its error bounds.

\subsection{An asymptotic consistent formulation}
Introduce
\be\label{Nd}
F^\ep(\bx,t)=N^\ep(\bx,t)+|E^\ep(\bx,t)|^2-G^\ep(\bx,t/\ep), \qquad \bx\in \mathbb{R}^d,\quad t\ge0,
\ee
where
\be\label{Gd}
G^\ep(\bx,s)=\ep^\alpha G_1(\bx,s)+\ep^{1+\beta}G_2(\bx,s),\qquad \bx\in \mathbb{R}^d,\quad s\ge0,
\ee
with $G_j(\bx,s)$ ($j=1,2$) being the solutions of the linear wave equations
\be\label{wave}
\begin{split}
&\p_{ss}G_j(\bx,s)-\Delta G_j(\bx,s)=0, \qquad \bx\in \mathbb{R}^d,\quad s>0,\\
&G_1(\bx,0)=\og_0(\bx), \quad \p_s G_1(\bx,0)\equiv 0, \quad
G_2(\bx,0)\equiv0,\quad  \p_s G_2(\bx,0)=\og_1(\bx).
\end{split}
\ee
Plugging (\ref{Nd}) into the ZS \eqref{Zak}, we can reformulate it into an
asymptotic consistent formulation
\be\label{Zak2}
\begin{split}
&i\p_tE^\ep(\bx,t)+\Dt E^\ep(\bx,t)+\left[|E^\ep(\bx,t)|^2-F^\ep(\bx,t)-
G^\ep(\bx,t/\ep)\right]E^\ep(\bx,t)=0,  \\
&\ep^2\p_{tt} F^\ep(\bx,t)-\Dt F^\ep(\bx,t)-\ep^2\p_{tt} |E^\ep(\bx,t)|^2=0,
\quad\bx \in \mathbb{R}^d, \quad t>0,\\
&E^\ep(\bx,0)=E_0(\bx),\quad F^\ep(\bx,0)\equiv0,\quad \p_t F^\ep (\bx,0)\equiv
0, \quad\bx\in \mathbb{R}^d.
\end{split}
\ee
Now the initial conditions in (\ref{Zak2}) are always well-prepared for any $\alpha\ge0$ and
$\beta\ge0$. In addition, the above system conserves the wave energy (\ref{mass})
and the `modified' {\sl Hamiltonian}
\bea
\wt{\mathcal L}^\ep(t)&:=&\int_{\mathbb{R}^d}\biggl[|\nabla E^\ep|^2
-\fl{1}{2}|E^\ep|^4+\frac{1}{2}|F^\ep|^2+\frac{1}{\ep^2}\int_0^t\int_0^s
\nabla F^\ep(\bx,s)\cdot \nabla F^\ep(\bx,s^\prime)ds^\prime ds\nonumber\\
&&+\int_0^t \left[G(\bx,s/\ep)\p_s |E^\ep(\bx,s)|^2
- \phi_1(\bx)F^\ep(\bx,s)\right]ds\biggr]d\bx\equiv \wt{\mathcal L}^\ep(0), \quad t\ge0.
\eea

When $\ep\to0^+$, i.e., in the subsonic limit regime, formally
we get $E^\ep(\bx,t)\to E(\bx,t)$ and $F^\ep(\bx,t)\to 0$, where $E:=E(\bx,t)$ satisfies
the NLSE \eqref{NLS}. In addition, when $\ep\to0^+$, formally we can also get
$E^\ep(\bx,t)\to \wt{E}^\ep(\bx,t)$ and $F^\ep(\bx,t)\to 0$, where
$\wt{E}^\ep:=\wt{E}^\ep(\bx,t)$ satisfies the following
nonlinear Schr\"{o}dinger equation with an oscillatory
potential $G^\ep(\bx,t/\ep)$ (NLSE-OP)
\be\label{ANLS}
\begin{split}
&i\p_t\wt{E}^\ep(\bx,t)+\Delta \wt{E}^\ep(\bx,t)+\left[|\wt{E}^\ep(\bx,t)|^2-
G^\ep(\bx,t/\ep)\right]
\wt{E}^\ep(\bx,t)=0,  \qquad t>0,\\
&\wt{E}^\ep(\bx,0)=E_0(\bx),\quad \bx\in \mathbb{R}^d.
\end{split}
\ee
It conserves the wave energy (\ref{mass}) with $E^\ep=\wt{E}^\ep$
and the `modified' {\sl Hamiltonian}
\bea\label{energy-nlsop}
\qquad \ \wt{{\mathcal L}}(t):=\int_{\mathbb{R}^d}\left[|\nabla \wt{E}^\ep|^2-\fl{1}{2}|\wt{E}^\ep|^4
+\int_0^t G(\bx,\fl{s}{\ep})\partial_s|\wt{E}^\ep(\bx,s)|^2ds
\right]d\bx
\equiv\wt{{\mathcal L}}(0), \ t\ge0.
\eea

\subsection{A uniformly accurate finite difference method}
For simplicity of notations, we will only present the numerical method
for the ZS \eqref{Zak2} in 1D and extensions to higher dimensions are straightforward.
When $d=1$, we truncate ZS on a bounded computational interval $\Og=(a,b)$
with homogeneous Dirichlet boundary condition (here
$|a|$ and $b$ are chosen large enough such that the truncation error
is negligible):
\be\label{Zak1d}
\begin{split}
&i\p_tE^\ep(x,t)+\p_{xx} E^\ep(x,t)+[|E^\ep(x,t)|^2-F^\ep(x,t)- G^\ep(x,t/\ep)]E^\ep(x,t)=0,  \\
&\ep^2\p_{tt} F^\ep(x,t)-\p_{xx} F^\ep(x,t)-\ep^2\p_{tt} |E^\ep(x,t)|^2=0,
\quad x \in \Omega, \quad t>0,\\
&E^\ep(x,0)=E_0(x),\quad F^\ep(x,0)\equiv0,\quad \p_t F^\ep (x,0)\equiv0, \quad x\in \bar\Omega,\\
&E^\ep(a,t)=E^\ep(b,t)=0,\qquad F^\ep(a,t)=F^\ep(b,t)=0, \qquad t\ge0,
\end{split}
\ee
where $G^\ep(x,s)$ is defined as \eqref{Gd} with $d=1$ and $G_j(x,s)$ ($j=1,2$)
being the solutions of the wave equations
\be\label{wave1D}
\begin{split}
&\p_{ss}G_j(x,s)-\partial_{xx} G_j(x,s)=0, \qquad x\in \Omega,\quad s>0,\\
&G_1(x,0)=\og_0(x), \quad \p_s G_1(x,0)\equiv 0, \quad
G_2(x,0)\equiv0,\quad  \p_s G_2(x,0)=\og_1(x),\\
&G_1(a,s)=G_1(b,s)=G_2(a,s)=G_2(b,s)=0, \qquad s\ge0.
\end{split}
\ee
When $\ep\to0^+$,  formally we get
$E^\ep(x,t)\to \wt{E}^\ep(x,t)$ and $F^\ep(x,t)\to 0$, where
$\wt{E}^\ep:=\wt{E}^\ep(x,t)$ satisfies the NLSE-OP
\bea\label{ANLS1D}
\begin{split}
&\qquad i\p_t\wt{E}^\ep(x,t)+\partial_{xx} \wt{E}^\ep(x,t)+\left[|\wt{E}^\ep(x,t)|^2-
G^\ep(x,t/\ep)\right]
\wt{E}^\ep(x,t)=0,\ \ t>0,\\
&\qquad \wt{E}^\ep(x,0)=E_0(x),\quad \bx\in \bar{\Omega}; \qquad \wt{E}^\ep(a,t)=\wt{E}^\ep(b,t)=0,\quad t\ge0.
\end{split}
\eea


Choose a mesh size $h:=\Dt x=(b-a)/M$ with $M$ being a positive integer and
a time step $\tau:=\Dt t>0$ and denote the grid points and time steps as
$$x_j:=a+jh,\quad j=0,1,\cdots,M;\quad t_k:=k\tau,\quad k=0,1,2,\cdots.$$
Define the index sets
$$\mathcal {T}_M=\{j \ | \ j=1,2,\cdots,M-1\},\quad \mathcal{T}_M^0=\{j\ |\ j=0,1,\cdots,M\}.$$
Let $E^{\ep,k}_j$ and $F^{\ep,k}_j$ be the approximations of $E^{\ep}(x_j,t_k)$ and
$F^\ep(x_j,t_k)$, respectively, and denote $E^{\ep,k}=(E^{\ep,k}_0,E^{\ep,k}_1,\ldots,E^{\ep,k}_M)^T\in \mathbb{C}^{(M+1)}$,
$F^{\ep,k}=(F^{\ep,k}_0,F^{\ep,k}_1,\ldots,F^{\ep,k}_M)^T\in \mathbb{R}^{(M+1)}$ as
the numerical solution vectors at $t=t_k$. Define the standard finite difference operators
\beas
&&\dt_t^+E_j^k=\fl{E_j^{k+1}-E_j^k}{\tau},\quad
\dt_t E_j^k=\fl{E_j^{k+1}-E_j^{k-1}}{2\tau},\quad \dt_t^2E_j^k=\fl{E_j^{k+1}-2E_j^k+E_j^{k-1}}{\tau^2},\\
&&\dt_x^+E_j^k=\fl{E_{j+1}^k-E_j^k}{h},\qquad
\dt_x^2E_j^k=\fl{E_{j+1}^k-2E_j^k+E_{j-1}^k}{h^2}.
\eeas

We present a finite difference discretization of \eqref{Zak1d} as following
\be\label{scheme1}
\begin{aligned}
&i\dt_t E_j^{\ep,k}=\left(-\dt_x^2-|E_j^{\ep,k}|^2+H_j^{\ep,k}
+\fl{F_j^{\ep,k+1}+F_j^{\ep,k-1}}{2}\right)\fl{E_j^{\ep,k+1}+E_j^{\ep,k-1}}{2},\\
&\ep^2\dt_t^2F_j^{\ep,k}=\fl{1}{2}\dt_x^2(F_j^{\ep,k+1}+F_j^{\ep,k-1})+
\ep^2\dt_t^2|E_j^{\ep,k}|^2, \quad j\in \mathcal{T}_M, \quad k\ge 1,
\end{aligned}
\ee
where an average of the oscillatory potential $G^\ep$ over the interval
$[t_{k-1},t_{k+1}]$ is used
\bea\label{Hd}
\qquad \ \ H_j^{\ep,k}=\fl{1}{2\tau}\int_{t_{k-1}}^{t_{k+1}}G^\ep(x_j,s/\ep)ds
=\fl{\ep}{2\tau}\int_{t_{k-1}/\ep}^{t_{k+1}/\ep}G^\ep(x_j,u)du, \quad j\in\mathcal{T}_M,\quad
 k\ge1.
\eea
The boundary and initial conditions are discretized as
\bea\label{initfd11}
\qquad \ \ E_0^{\ep,k}=E_M^{\ep,k}=F_0^{\ep,k}=F_M^{\ep,k}=0,\ k\ge0; \quad
E_j^{\ep,0}=E_0(x_j),\ F_j^{\ep,0}=0, \quad j\in \mathcal{T}_M^0.
\eea
In addition, the first step $E_j^{\ep,1}$ and $F_j^{\ep,1}$ can be obtained via (\ref{Zak1d}) and
the Taylor expansion as
\be\label{F_1}
E_j^{\ep,1}= E_0(x_j)+\tau \phi_2(x_j)+\fl{\tau^2}{2}\phi_3(x_j),\quad
F_j^{\ep,1}=\fl{\tau^2}{2}\phi_4(x_j),\quad j\in \mathcal{T}_M,
\ee
where
\be\label{initFD1d}
\begin{split}
&\phi_2(x):=\p_t E^\ep(x,0)=i\left[E_0^{\prime\prime}(x)-N_0^\ep(x)E_0(x)\right],\\
&\phi_3(x):=\p_{tt}E^\ep(x,0)=i\left[\phi_2^{\prime\prime}(x)-N_1^\ep(x)E_0(x)-N_0^\ep(x)\phi_2(x)\right],\qquad x\in\Omega,\\
&\phi_4(x):=\p_{tt}F^\ep(x,0)=\p_{tt}\rho^\ep(x,0)=2{\rm Im}\left[\phi_2(x)\overline{E_0^{\prime\prime}(x)}
+E_0(x)\overline{\phi_2^{\prime\prime}(x)}\right].
\end{split}
\ee

   If it is needed in practical computation, the second order derivatives in (\ref{initFD1d})
can be approximated by the second order finite difference as
$f^{\prime\prime}(x_j)\approx \delta_x^2 f(x_j)$ for $j\in\mathcal{T}_M$.
In addition, $H_j^{\ep,k}$ in (\ref{Hd}) can be approximated by solving the wave
equations (\ref{wave1D}) via the sine pseudospectral method in space
and then integrating in time in phase space {\sl exactly} as
\beas
H_j^{\ep,k}&\approx&\frac{\ep}{2\tau}\sum_{l=1}^{M-1}\sin\left(\mu_l(x_j-a)\right)
\int_{t_{k-1}/\ep}^{t_{k+1}/\ep}
\left[\ep^\alpha\wt{(\omega_0)}_l\cos(\mu_l u) +\frac{\ep^{1+\beta}}{\mu_l}
\wt{(\omega_1)}_l\sin(\mu_l u)\right]du\nonumber\\
&=&\sum_{l=1}^{M-1}\frac{\ep }{\tau \mu_l}\sin\left(\frac{lj\pi}{M}\right)\sin\left(\frac{\tau\mu_l}{\ep}\right)
\left[\ep^\alpha \wt{(\omega_0)}_l\cos\left(\frac{\mu_l t_{k}}{\ep}\right)+\frac{\ep^{1+\beta}}{\mu_l}\wt{(\omega_1)}_l
\sin\left(\frac{\mu_l t_{k}}{\ep}\right)\right],
\eeas
where for $l=1,2,\ldots,M-1$,
\[
\mu_l=\frac{l\pi}{b-a}, \quad \wt{(\omega_0)}_l=\frac{2}{M}\sum_{j=1}^{M-1}\omega_0(x_j)
\sin\left(\frac{lj\pi}{M}\right), \quad \wt{(\omega_1)}_l=\frac{2}{M}\sum_{j=1}^{M-1}\omega_1(x_j)
\sin\left(\frac{lj\pi}{M}\right).
\]

\subsection{Main results}
For convenience of notation, denote
$$0\le \alpha^*=\min \{\alpha,1\}\le 1.$$
Let $T^*>0$ be the maximum common existence time for the solutions of the ZS \eqref{Zak1d} and
the NLSE-OP \eqref{ANLS1D}. Then
for any fixed $0<T<T^*$, according to the known results
in \cite{Added, Masmoudi2008,Ozawa, Schochet}, we assume that
the solution $(E^\ep, F^\ep)$ of the ZS \eqref{Zak1d} and the solution
$\wt{E}^\ep$ of the NLSE-OP \eqref{ANLS1D} are smooth enough over $\Omega_T:=\Omega\times [0,T]$ and satisfy
\begin{equation*}
(A)\begin{split}
&\|E^\ep\|_{W^{5,\infty}}+\|\p_t E^\ep\|_{W^{1,\infty}}
+\|\p_{tt} F^\ep\|_{W^{2,\infty}}
+\|\wt{E}^\ep\|_{W^{5,\infty}} +\|\p_t\wt{E}^\ep\|_{W^{1,\infty}}\lesssim 1,\\
&\|F^\ep\|_{W^{4,\infty}}\lesssim \ep^2,\qquad
\|\p_t F^\ep\|_{W^{4,\infty}}\lesssim \ep, \qquad
\|\p_{tt}\wt{E}^\ep\|_{W^{4,\infty}}\lesssim \frac{1}{\ep^{1-\alpha^*}},\\
&\|\p_{tt}E^\ep\|_{W^{4,\infty}}+\|\p_t^3 F^\ep\|_{W^{2,\infty}}\lesssim \frac{1}{\ep},
\qquad\|\p_t^3 E^\ep\|_{W^{4,\infty}}+\|\p_t^4 F^\ep\|_{W^{2,\infty}}\lesssim \frac{1}{\ep^2}.
\end{split}
\end{equation*}
We further assume that  the initial data satisfy
\begin{equation*}
(B)\hskip3cm
\|E_0\|_{W^{5,\infty}(\Omega)}+\|\og_0\|_{W^{3,\infty}(\Omega)}+
\|\og_1\|_{W^{3,\infty}(\Omega)}\lesssim 1.\hskip6cm
\end{equation*}
Then one can obtain \cite{Masmoudi2008,Schochet,Sulemb}
\be\label{Gp}
\|G^\ep\|_{W^{3,\infty}([0,\infty),W^{3,\infty}(\Omega))}\lesssim
\ep^{\alpha^*}.
\ee
In addition, we assume the following convergence rate from ZS to NLSE-OP
\begin{equation*}
(C)\hskip3cm
\|E^\ep-\wt{E}^\ep\|_{L^\infty([0,T];H^1(\Omega))}\lesssim \ep^2.\hskip6cm
\end{equation*}

 Denote
$$X_M=\left\{v=\left(v_0,v_1,\ldots,v_M\right)^T\  | \ v_0=v_M=0\right\} \subseteq \mathbb{C}^{M+1},$$
equipped with norms and inner products defined as
\beas
&&\|u\|^2=h\sum\limits_{j=1}^{M-1}|u_j|^2, \quad\|\dt_x^+ u\|^2=h\sum\limits_{j=0}^{M-1} |\dt_x^+ u_j|^2, \quad \|u\|_\infty=\sup\limits_{j\in \mathcal{T}_M^0}|u_j|,\\
&&(u,v)=h\sum\limits_{j=1}^{M-1}u_j \overline{v_j}, \quad
\langle \delta_x^+u,\delta_x^+v\rangle=h\sum\limits_{j=0}^{M-1} (\delta_x^+u_j) \;
(\delta_x^+\overline{v_j}), \qquad u,v\in X_M.
\eeas
Then we have
\be\label{innpX_M}
(-\dt_x^2 u,v)=\langle\dt_x^+ u,\dt_x^+ v\rangle,\quad
((-\dt_x^2)^{-1}u,v)=(u,(-\dt_x^2)^{-1}v), \qquad u, v\in X_M.
\ee

Define the error functions $e^{\ep,k}\in X_M$ and $f^{\ep,k}\in X_M$ as
\be
e^{\ep,k}_j=E^\ep(x_j,t_k)-E_j^{\ep,k},\quad f_j^{\ep,k}=F^\ep(x_j,t_k)
-F_j^{\ep,k},\quad j\in \mathcal{T}_M^0,\quad 0\le k\le \frac{T}{\tau}.
\ee
Then we have the following error estimates for (\ref{scheme1}) with (\ref{Hd})-(\ref{F_1}).
\begin{theorem}\label{thm1}
Under the assumptions (A)-(C), there exist $h_0>0$ and $\tau_0>0$ sufficiently small
and independent of $0<\ep\le 1$ such that, when $0<h\le h_0$ and $0<\tau\le \tau_0$,
the following two error estimates of the scheme \eqref{scheme1} with \eqref{Hd}-\eqref{F_1}
hold
\begin{align}
\|e^{\ep,k}\|+\|\dt_x^+e^{\ep,k}\|+\|f^{\ep,k}\|&\lesssim h^2+ \fl{\tau^2}{\ep},
 \qquad 0\le k\le \fl{T}{\tau},\quad  0<\ep\le 1,\label{esti1}\\
\|e^{\ep,k}\|+\|\dt_x^+e^{\ep,k}\|+\|f^{\ep,k}\|&\lesssim h^2+\tau^2+\tau\ep^{\alpha^*}
+\ep^{1+\alpha^*}.\label{esti2}
\end{align}
Thus by taking the minimum among the two error bounds for $\ep\in (0,1]$, we
obtain a uniform error estimate for well-prepared initial data, i.e., $\alpha\ge1$,
\be\label{uniform_wp}
\|e^{\ep,k}\|+\|\dt_x^+e^{\ep,k}\|+\|f^{\ep,k}\|\lesssim h^2+\min_{0<\ep\le1}\left\{\tau^2+\tau\ep+\ep^2,\frac{\tau^2}{\ep}\right\}\lesssim h^2+\tau^{4/3},
\ee
and respectively, for ill-prepared initial data, i.e., $0\le \alpha<1$,
\be\label{uniform_ip}
\|e^{\ep,k}\|+\|\dt_x^+e^{\ep,k}\|+\|f^{\ep,k}\|\lesssim h^2+\min_{0<\ep\le1}\left\{\tau^2+\ep^\alpha(\tau+\ep),
\frac{\tau^2}{\ep}\right\}\lesssim h^2+\tau^{1+\frac{\alpha}{2+\alpha}}.
\ee
\end{theorem}

\section{Error analysis}
In order to prove Theorem \ref{thm1}, we will use the energy method to obtain one error bound
\eqref{esti1}
and use the limiting equation NLSE-OP (\ref{ANLS1D}) to get the other one \eqref{esti2},
which is shown in the following diagram \cite{Cai2012,Cai2013,Bao2014,Jin1999,Degond}.
\[
\xymatrixcolsep{8pc}\xymatrix{
(E^{\ep,k}, F^{\ep,k}) \ar[r]^{\quad O(h^2+\tau^2+
\tau\ep^{\alpha^*}+\ep^{1+\alpha^*})} \ar[rd]^{ 
}_{O(h^2+\tau^2/\ep)}&
(\wt{E}^\ep,0)\ar[d]^{O(\ep^2)} \\
&(E^\ep,F^\ep)}
\]

To simplify notations, for a function $V:=V(x,t)$ and a grid function $V^{k}\in X_M$ with $k\ge0$,
we denote for $k\ge1$
\[\Lparen V\Rparen(x,t_k)=\fl{V(x,t_{k+1})+V(x,t_{k-1})}{2},\quad x\in\bar{\Omega}; \quad
\Lbrack V\Rbrack_j^k=\fl{V^{k+1}_j+V^{k-1}_j}{2}, \quad j\in \mathcal{T}_M^0.
\]

In order to deal with the nonlinearity and to bound the numerical solution,
we adapt the cut-off technique which has been widely used in the literatures
\cite{Akrivis, Cai2013, Bao2012, Thomee}, i.e.
the nonlinearity is first truncated to a global Lipschitz function with
compact support and then the error bound can be achieved if the
exact solution is bounded and the numerical solution is close to
the exact solution under some conditions on the mesh size and time step.
Choose  a smooth function $\gamma(s)\in C^\infty(\mathbb{R})$ such that
$$\gamma(s)=
\left\{
\begin{aligned}
&1, \quad &|s|\le 1,\\
&\in[0,1],\quad &|s|\le 2,\\
&0,\quad &|s|\ge 2, \\
\end{aligned}
\right.
$$
and by assumption (A) we can choose  $M_0>0$ as
$$M_0=\max \left\{\sup\limits_{\ep\in(0,1]}
\|E^\ep\|_{L^\infty(\Omega_T)}, \sup\limits_{\ep\in(0,1]}
\|\wt{E}^\ep\|_{L^\infty(\Omega_T)}\right\}. $$
For $s\ge 0$, $y_1, y_2\in \mathbb{C}$, define
$$\gamma_{_B}(s)=s\,\gamma\left(\frac{s}{B}\right),\quad \hbox{with}\quad  B=(M_0+1)^2,$$
and
$$g(y_1,y_2)=\fl{y_1+y_2}{2}\int_0^1\gamma_{_B}'(s|y_1|^2+(1-s)|y_2|^2)ds
=\fl{\gamma_{_B}(|y_1|^2)-\gamma_{_B}(|y_2|^2)}{|y_1|^2-|y_2|^2}\cdot\fl{y_1+y_2}{2}.
$$
Then $\gamma_{_B}(s)$ is global Lipschitz and there exists $C_{_B}>0$, such that
\be\label{rhoB}
|\gamma_{_B}(s_1)-\gamma_{_B}(s_2)|\le \sqrt{C_{_B}}|\sqrt{s_1}-\sqrt{s_2}|,\quad \forall s_1, s_2 \ge 0.
\ee

Let $\hat{E}^{\ep,k}$, $\hat{F}^{\ep,k}\in X_M$ ($k\ge 0$) be the solution of the following
\be\label{scheme}
\begin{aligned}
&i\dt_t \hat{E}_j^{\ep,k}=(-\dt_x^2+H_j^{\ep,k})
\Lbrack\hat{E}^{\ep}\Rbrack_j^k+\big(-\gamma_{_B}(|\hat{E}_j^{\ep,k}|^2)
+\Lbrack\hat{F}^\ep\Rbrack_j^k\big)g(\hat{E}^{\ep,k+1}_j,\hat{E}_j^{\ep,k-1}),\\
&\ep^2\dt_t^2\hat{F}_j^{\ep,k}=\fl{1}{2}\dt_x^2(\hat{F}_j^{\ep,k+1}+\hat{F}_j^{\ep,k-1})+
\ep^2\dt_t^2\gamma_{_B}(|\hat{E}_j^{\ep,k}|^2), \qquad j\in \mathcal{T}_M, \quad k\ge1,\\
&\hat{E}^{\ep,0}_j=E^{\ep,0}_j, \quad \hat{F}^{\ep,0}_j=F^{\ep,0}_j=0,\quad
\hat{E}^{\ep,1}_j=E^{\ep,1}_j, \quad \hat{F}^{\ep,1}_j=F^{\ep,1}_j,\qquad j\in \mathcal{T}_M^0.
\end{aligned}
\ee
Here $(\hat{E}^{\ep,k},\hat{F}^{\ep,k})$ can be viewed as another approximation of
the solution  $(E^\ep,F^\ep)$ of ZS with a cut-off Lipschitz nonlinearity.
Define error functions $\hat{e}^{\ep,k}$, $\hat{f}^{\ep,k}\in X_M$ as
\be
\hat{e}_j^{\ep,k}=E^\ep(x_j,t_k)-\hat{E}^{\ep,k}_j,\quad
\hat{f}_j^{\ep,k}=F^{\ep}(x_j,t_k)-\hat{F}^{\ep,k}_j, \quad j\in \mathcal{T}_M^0, \quad k\ge0.
\ee

For $(\hat{e}^{\ep,k},\hat{f}^{\ep,k})$, we have the following estimates.
\begin{theorem}\label{thm2}
Under the assumption (A), there exists $\tau_1>0$
sufficiently small and independent of $0<\ep\le1$ such that, when $0<\tau\le \tau_1$ and
$0<h\le \frac{1}{2}$, we have the following error estimate for
the scheme \eqref{scheme}
\be \label{errort1}
\|\hat{e}^{\ep,k}\|+\|\dt_x^+\hat{e}^{\ep,k}\|+\|\hat{f}^{\ep,k}\|\lesssim h^2+
\fl{\tau^2}{\ep}, \quad 0\le k\le \fl{T}{\tau}, \quad 0<\ep\le1.
\ee
\end{theorem}

Introduce local truncation errors $\hat{\xi}_j^{\ep,k}$, $\hat{\eta}_j^{\ep,k}\in X_M$ as
\be\label{local}
\begin{split}
\hat{\xi}_j^{\ep,k}=&i\dt_t E^\ep(x_j,t_k)+(\dt_x^2-H_j^{\ep,k})\Lparen E^\ep\Rparen(x_j,t_k)\\
&+\big(\gamma_{_B}(|E^\ep(x_j,t_k)|^2)-\Lparen F^\ep\Rparen(x_j,t_k)\big)
g\big(E^\ep(x_j,t_{k+1}),E^\ep(x_j,t_{k-1})\big)\\
=&i\dt_t E^\ep(x_j,t_k)+(\dt_x^2+|E^\ep(x_j,t_{k})|^2-H_j^{\ep,k}-
\Lparen F^\ep\Rparen(x_j,t_k))\Lparen E^\ep\Rparen(x_j,t_k),\\
\hat{\eta}_j^{\ep,k}=&\ep^2\dt_t^2F^\ep(x_j,t_k)-\dt_x^2\Lparen F^\ep\Rparen(x_j,t_k)-
\ep^2\dt_t^2\gamma_{_B}(|E^\ep(x_j,t_k)|^2)\\
=&\ep^2\dt_t^2F^\ep(x_j,t_k)-\dt_x^2\Lparen F^\ep\Rparen(x_j,t_k)-
\ep^2\dt_t^2|E^\ep(x_j,t_k)|^2, \quad j \in \mathcal{T}_M, \quad k\ge1.
\end{split}
\ee
Then we have
\begin{lemma}\label{local_e}
Under the assumption (A), when $0<h\le \frac{1}{2}$ and $0<\tau\le \frac{1}{2}$,  we have
\be
|\hat{\xi}_j^{\ep,k}|+|\dt_x^+\hat{\xi}_j^{\ep,k}|\lesssim h^2+\fl{\tau^2}{\ep},\ \
|\hat{\eta}_j^{\ep,k}|\lesssim \ep^2h^2+\tau^2,\ \
|\dt_t\hat{\eta}_j^{\ep,k}|\lesssim \ep h^2 +\fl{\tau^2}{\ep},\quad  j \in \mathcal{T}_M.
\ee
\end{lemma}
\emph{Proof.} By \eqref{Zak1d} and using Taylor expansion, we get
\begin{align*}
&i\dt_t E^\ep(x_j,t_k)=\fl{i}{2\tau}\int_{t_{k-1}}^{t_{k+1}}\p_t E^\ep(x_j,s)ds\\
&=\fl{1}{2\tau}\int_{t_{k-1}}^{t_{k+1}}\left[(-\p_{xx}E^\ep-|E^\ep|^2E^\ep
+ E^\ep F^\ep)(x_j,s)+E^\ep(x_j,s) G^\ep\left(x_j,\fl{s}{\ep}\right)\right]ds\\
&=-E_{xx}^\ep(x_j,t_k)-|E^\ep(x_j,t_k)|^2E^\ep(x_j,t_k)+ E^\ep(x_j,t_k)F^\ep(x_j,t_k)\\
&\quad-\fl{\tau^2}{4}\int_{0}^1(1-s)^2\sum\limits_{m=\pm 1}\p_{tt}
(E_{xx}^\ep+|E^\ep|^2E^\ep-E^\ep F^\ep)(x_j,t_k+ms \tau)ds\\
&\quad+\fl{1}{2\tau}\int_{-\tau}^\tau E^\ep(x_j,t_k+s) G^\ep\left(x_j,
\fl{t_k+s}{\ep}\right)ds, \quad j \in \mathcal{T}_M, \quad 1\le k\le \frac{T}{\tau}-1.
\end{align*}
Similarly, by Taylor expansion, we have
\begin{align*}
&\left(\dt_x^2+|E^\ep(x_j,t_{k})|^2-H_j^{\ep,k}-\Lparen F^\ep\Rparen(x_j,t_k)\right)\Lparen E^\ep\Rparen(x_j,t_k)\\
&=E^\ep_{xx}(x_j,t_k)+\big(|E^\ep(x_j,t_{k})|^2-H_j^{\ep,k}-F^\ep(x_j,t_k)\big)
E^\ep(x_j,t_k)\\
&\ \ +\fl{h^2}{12}\int_0^1(1-s)^3\sum\limits_{m=\pm1}\sum\limits_{l=\pm1}
E^\ep_{xxxx}(x_j+slh,t_k+m \tau)ds\\
&\ \ +\fl{\tau^2}{2}\int_0^1(1-s)\sum\limits_{m=\pm1}
\Bigl(E_{xxtt}^\ep(x_j,t_k+ms\tau)-E^\ep(x_j,t_k)F_{tt}^\ep(x_j,t_k+ms\tau)\Bigr)
ds\\
&\ \ +\fl{\tau^2}{2}\left(|E^\ep(x_j,t_{k})|^2-H_j^{\ep,k}-\Lparen F^\ep \Rparen(x_j,t_k)\right)
\int_0^1(1-s)\sum\limits_{m=\pm1}E^\ep_{tt}(x_j,t_k+ms\tau)ds.
\end{align*}
Note that by \eqref{Hd}, we have
\begin{align*}
&\fl{1}{2\tau}\int_{-\tau}^\tau E^\ep(x_j,t_k+s) G^\ep\left(x_j,\fl{t_k+s}{\ep}\right)ds
-E^\ep(x_j,t_k)H_j^{\ep,k}\\
&=\fl{1}{2\tau}E^\ep_t(x_j,t_k)\int_{-\tau}^\tau s\, G^\ep\left(x_j,\fl{t_k+s}{\ep}\right)ds+A_1\\
&=\fl{\tau^2}{2}E^\ep_t(x_j,t_k)\int_0^1\fl{s}{\ep}\int_{-s}^s
G_t^\ep\left(x_j,\fl{t_k+\tht \tau}{\ep}\right)d\tht ds+A_1,
\end{align*}
where
$$A_1=\fl{\tau^2}{2}\int_{-1}^1\int_0^s (s-\tht)\,G^\ep\left(x_j,\fl{t_k+s\tau}{\ep}\right)
E_{tt}^\ep(x_j,t_k+\tht\tau)
d\tht ds.$$
Accordingly, by the assumption (A) and \eqref{Gp}, we conclude that
\begin{align*}
|\hat{\xi}_j^{\ep,k}|&\lesssim h^2\|E^\ep_{xxxx}\|_{L^\infty}
+\tau^2
\Bigl[
\|E^\ep_{xxtt}\|_{L^\infty}+\|\p_{tt}(|E^\ep|^2E^\ep)\|_{L^\infty}
+\|E^\ep\|_{L^\infty}\|F_{tt}^\ep\|_{L^\infty}\\
&+\fl{1}{\ep}\|E_t^\ep\|_{L^\infty}\left(\|G_t^\ep\|_{L^\infty}+\ep
\|F_{t}^\ep\|_{L^\infty}\right)
+|E_{tt}^\ep\|_{L^\infty}\left(\|G^\ep\|_{L^\infty}+\|F^\ep\|_{L^\infty}+
\|E^\ep\|^2_{L^\infty}\right)\Bigr]\\
&\lesssim h^2+\fl{\tau^2}{\ep}, \qquad j \in \mathcal{T}_M, \quad 1\le k\le \frac{T}{\tau}-1.
\end{align*}
Applying $\dt_x^+$ to $\hat{\xi}^{\ep,k}$ and using the same approach, we get
$$|\dt_x^+\hat{\xi}_j^{\ep,k}|
\lesssim h^2+\fl{\tau^2}{\ep}, \qquad j \in \mathcal{T}_M, \quad 1\le k\le \frac{T}{\tau}-1.$$
Similarly, we obtain
\begin{align*}
\hat{\eta}_j^{\ep,k}&=\fl{\ep^2\tau^2}{6}\int_0^1(1-s)^3\sum\limits_{m=\pm1}
\left(F_{tttt}^\ep(x_j,t_k+ms\tau)-(|E^\ep|^2)_{tttt}(x_j,t_k+ms\tau)\right)ds\\
&\quad-\fl{\tau^2}{2}\int_0^1(1-s)\sum\limits_{m=\pm1}
F_{xxtt}^\ep(x_j,t_k+ms\tau)ds\\
&\quad-\fl{h^2}{12}\int_0^1(1-s)^3\sum\limits_{l=\pm1}\sum\limits_{m=\pm1}
F_{xxxx}^\ep(x_j+lsh,t_k+m\tau)ds,
\end{align*}
which implies
\begin{align*}
|\hat{\eta}_j^{\ep,k}|&\lesssim h^2\|F^\ep_{xxxx}\|_{L^\infty}
+\tau^2(\|F_{xxtt}^\ep\|_{L^\infty}+\ep^2\|F_{tttt}^\ep\|_{L^\infty
}+\ep^2\|\p_{tttt}|E^\ep|^2\|_{L^\infty})\\
&\lesssim \ep^2h^2+\tau^2, \qquad j \in \mathcal{T}_M, \quad 1\le k\le \frac{T}{\tau}-1.
\end{align*}
Applying $\dt_t$ to $\hat{\eta}_j^{\ep,k}$, we have
\begin{align*}
|\dt_t\hat{\eta}_j^{\ep,k}|&\lesssim h^2\|F^\ep_{xxxxt}\|_{L^\infty}
+\tau^2(\|F_{xxttt}^\ep\|_{L^\infty}+\ep^2\|F_{ttttt}^\ep\|_{L^\infty}+
\ep^2\|\p_{ttttt}|E^\ep|^2\|_{L^\infty})\\
&\lesssim \ep h^2+\fl{\tau^2}{\ep}, \qquad j \in \mathcal{T}_M, \quad 2\le k\le \frac{T}{\tau}-2.
\end{align*}
Thus the proof is completed.
\hfill $\square$

For the initial step, we have the following estimates.
\begin{lemma} \label{initial-e}
Under the assumption (A), when $0<\tau\le\frac{1}{2}$,
the first step errors of the discretization \eqref{scheme}
with \eqref{initfd11} and \eqref{F_1} satisfy
\bea\label{initial_e}
\qquad \ \ \hat{e}_j^{\ep,0}=\hat{f}_j^{\ep,0}=0,\ \
|\hat{e}_j^{\ep,1}|+|\dt_t^+\hat{f}_j^{\ep,0}|+|\dt_x^+\hat{e}_j^{\ep,1}|\lesssim \fl{\tau^2}{\ep},\ \
|\hat{f}_j^{\ep,1}|\lesssim \fl{\tau^3}{\ep},\ \
|\dt_t^+\hat{e}_j^{\ep,0}|\lesssim \fl{\tau^2}{\ep^2}.
\eea
\end{lemma}
\emph{Proof.}
By the definition of $\hat{E}^{\ep,1}_j$, we obtain
$$|\hat{e}_j^{\ep,1}|=\tau^2\,\left|\int_0^1(1-s)E_{tt}^\ep(x_j,s\tau)ds-
\fl{1}{2}E^\ep_{tt}(x_j,0)\right|\lesssim\tau^2\|E^\ep_{tt}\|_{L^\infty}
  \lesssim \fl{\tau^2}{\ep}.$$
On the other hand, we also have
$$ |\hat{e}_j^{\ep,1}|=\fl{\tau^3}{2}\, \left|\int_0^1(1-s)^2E_{ttt}^\ep(x_j,s\tau)ds\right|\lesssim\tau^3\|E^\ep_{ttt}\|_{L^\infty}
  \lesssim \fl{\tau^3}{\ep^2},$$
  which implies that $|\dt_t^+\hat{e}_j^{\ep,0}|\lesssim \fl{\tau^2}{\ep^2}$.
Similarly, $|\dt_x^+\hat{e}_j^{\ep,1}|\lesssim \tau^2 \|E_{xtt}^\ep\|_{L^\infty}\lesssim \fl{\tau^2}{\ep}$.
It follows from \eqref{F_1} and assumption (A) that
$$|\hat{f}^{\ep,1}_j|=\fl{\tau^3}{2}\left|\int_0^1(1-s)^2F_{ttt}^\ep(x_j,s\tau)ds\right|\lesssim
\tau^3\|F_{ttt}^\ep\|_{L^\infty}
  \lesssim \fl{\tau^3}{\ep}.$$
Recalling that $\hat{f}_j^{\ep,0}=0$, we can get that
$|\dt_t^+\hat{f}_j^{\ep,0}|\lesssim \fl{\tau^2}{\ep}$, which completes the proof.
\hfill $\square$

Subtracting \eqref{scheme} from \eqref{local}, we have the error
equations
\be
\begin{split}
&i\dt_t\hat{e}^{\ep,k}_j=(-\dt_x^2+H_j^{\ep,k})
\fl{\hat{e}_j^{\ep,k+1}+\hat{e}_j^{\ep,k-1}}{2}+r_j^k+\hat{\xi}_j^{\ep,k},\\
&\ep^2\dt_t^2\hat{f}_j^{\ep,k}=\fl{1}{2}\dt_x^2(\hat{f}_j^{\ep,k+1}+
\hat{f}_j^{\ep,k-1})+\ep^2\dt_t^2p_j^k+\hat{\eta}_j^{\ep,k},\quad
j \in \mathcal{T}_M, \ 1\le k\le \frac{T}{\tau}-1,
\end{split}
\label{n_eq}
\ee
where $r^k\in X_M$ and $p^k\in X_M$ are defined as
\be\label{pjk234}
\begin{split}
&r_j^k=\left.\left[-|E^\ep|^2
+\Lparen F^\ep\Rparen\right]\Lparen E^\ep\Rparen\right|_{(x_j,t_k)}
+\left[\gamma_{_B}(|\hat{E}_j^{\ep,k}|^2)
-\Lbrack\hat{F}^\ep\Rbrack_j^k\right] g(\hat{E}_j^{\ep,k+1},\hat{E}_j^{\ep,k-1}),\\
&p_j^k=|E^\ep(x_j,t_k)|^2-\gamma_{_B}(|\hat{E}_j^{\ep,k}|^2),\qquad
j \in \mathcal{T}_M, \quad 1\le k\le \frac{T}{\tau}-1.
\end{split}
\ee
By the property of $\gamma_{_B}$ in \eqref{rhoB}, we get for $0\le k\le \frac{T}{\tau}$
\be\label{p_k}
|p_j^k|=\left|\gamma_{_B}(|E^\ep(x_j,t_k)|^2)-\gamma_{_B}(|\hat{E}_j^{\ep,k}|^2)\right|\le \sqrt{C_{_B}} |\hat{e}_j^{\ep,k}|, \quad j \in \mathcal{T}_M.
\ee
Recalling the definition of $g(\cdot,\cdot)$ and noting
that $\Lparen E^\ep\Rparen(x_j,t_k)=g\left(E^\ep(x_j,t_{k+1}),
E^\ep(x_j,t_{k-1})\right)$, similar to the proof in \cite{Cai2012, Cai2015} with the details omitted
here for brevity, we have for $j\in \mathcal{T}_M$ and $1\le k\le \frac{T}{\tau}-1$,
\be\label{nonlbd}
\begin{split}
&\left|g(\hat{E}^{\ep,k+1}_j,\hat{E}_j^{\ep,k-1})\right|\lesssim 1, \qquad  \left|\Lparen E^\ep\Rparen(x_j,t_k)-g(\hat{E}_j^{\ep,k+1},\hat{E}_j^{\ep,k-1})\right|
\lesssim \sum_{l=k\pm1} |\hat{e}_j^{\ep,l}|,\\
&\left|\dt_x^+\big(\Lparen E^\ep\Rparen(x_j,t_k)
-g(\hat{E}_j^{\ep,k+1},\hat{E}_j^{\ep,k-1}))\right|\lesssim \sum\limits_{l=k\pm1}(|\hat{e}_j^{\ep,l}|+|\hat{e}_{j+1}^{\ep,l}|+
|\dt_x^+\hat{e}_j^{\ep,l}|).
\end{split}
\ee

\emph{Proof of Theorem \ref{thm2}.}
Multiplying both sides of the first equation in \eqref{n_eq} by $4\tau\,
\overline{\Lbrack\hat{e^\ep}\Rbrack_j^k}$, summing together for $j\in \mathcal{T}_M$ and
taking the imaginary parts, we obtain
\be\label{e_eq1}
\|\hat{e}^{\ep,k+1}\|^2-\|\hat{e}^{\ep,k-1}\|^2
=4\tau\,\mathrm{Im}\left(r^k+\hat{\xi}^{\ep,k},
\Lbrack\hat{e}^\ep\Rbrack^k\right), \qquad 1\le k\le \frac{T}{\tau}-1.
\ee
Using the same approach by multiplying $4\tau\dt_t\overline{\hat{e}_j^{\ep,k}}$ and taking the real parts,
 we get
\be\label{e_eq2}
\|\dt_x^+\hat{e}^{\ep,k+1}\|^2-\|\dt_x^+\hat{e}^{\ep,k-1}\|^2
=-4\,\mathrm{Re}\left(H^{\ep,k}\Lbrack\hat{e}^\ep\Rbrack^k+r^k+
\hat{\xi}^{\ep,k}, \tau\dt_t\hat{e}^{\ep,k}\right).
\ee
Introduce $\hat{u}^{\ep,k+1/2}\in X_M$ satisfying
$$-\dt_x^2\hat{u}_j^{\ep,k+1/2}=\dt_t^+(\hat{f}_j^{\ep,k}-p_j^k),\qquad j\in \mathcal{T}_M.$$
Multiplying both sides of the second equation in \eqref{n_eq} by $\tau(\hat{u}_j^{\ep,k+1/2}+
\hat{u}_j^{\ep,k-1/2})$,
summing together for $j\in \mathcal{T}_M$, we obtain
\bea
&&\ep^2\left(\|\dt_x^+\hat{u}^{\ep,k+1/2}\|^2-\|\dt_x^+\hat{u}^{\ep,k-1/2}\|^2\right)
+\fl{1}{2}\left(\|\hat{f}^{\ep,k+1}\|^2-\|\hat{f}^{\ep,k-1}\|^2\right)\nn\\
&&=\left(\Lbrack\hat{f}^\ep\Rbrack^k, 2\tau \dt_t p^{k}\right)
+\tau\left(\hat{\eta}^{\ep,k},\hat{u}^{\ep,k+1/2}+\hat{u}^{\ep,k-1/2}\right),
\qquad 1\le k\le \frac{T}{\tau}-1.\label{n_eq2}
\eea
Define a discrete `energy'
\bea
\mathcal {A}^k&=&C_{_B}(\|\hat{e}^{\ep,k}\|^2+\hat{e}^{\ep,k+1}\|^2)+
\|\dt_x^+\hat{e}^{\ep,k}\|^2+\|\dt_x^+\hat{e}^{\ep,k+1}\|^2\nonumber\\
&&+\ep^2\|\dt_x^+\hat{u}^{\ep,k+1/2}\|^2
+\fl{1}{2}(\|\hat{f}^{\ep,k}\|^2+\|\hat{f}^{\ep,k+1}\|^2),
\qquad 0\le k\le \frac{T}{\tau}-1.\label{energy_d}
\eea
Multiplying \eqref{e_eq1} by $C_{_B}>0$ and then summing with \eqref{e_eq2} and \eqref{n_eq2}, we get
\bea\label{sp}
\mathcal{A}^k-\mathcal{A}^{k-1}&=&4\tau C_{_B}\;\mathrm{Im}\left(r^k+\hat{\xi}^{\ep,k},
\Lbrack\hat{e}^\ep\Rbrack^k\right)-4\,\mathrm{Re}
\left(H^{\ep,k}\Lbrack\hat{e}^\ep\Rbrack^k+r^k+\hat{\xi}^{\ep,k},\tau \dt_t\hat{e}^{\ep,k}\right)\nonumber\\
&&+\left(\Lbrack\hat{f}^\ep\Rbrack^k, 2\tau \dt_t p^{k}\right)+\tau\left(\hat{\eta}^{\ep,k},\hat{u}^{\ep,k+1/2}+\hat{u}^{\ep,k-1/2}\right),
\quad 1\le k\le \frac{T}{\tau}-1.
\eea

Now we estimate different terms in the right hand side of \eqref{sp}.
Let $q_1^k\in X_M$ and  $q_2^k\in X_M$ defined as
\be\label{q1}
\begin{split}
&q^k_{1j}=\big(-|E^\ep(x_j,t_{k})|^2
+\Lparen F^\ep\Rparen(x_j,t_k)\big)\big(\Lparen E^\ep\Rparen(x_j,t_k)-
g(\hat{E}^{\ep,k+1}_j,\hat{E}_j^{\ep,k-1})\big),\\
&q_{2j}^k=-g(\hat{E}^{\ep,k+1}_j,\hat{E}_j^{\ep,k-1})\big(p_j^{k}
-\Lbrack\hat{f}^\ep\Rbrack_j^k\big), \qquad j\in \mathcal{T}_M.
\end{split}
\ee
Then we have
\be\label{r_d}
r^k=q_{1}^k+q_{2}^k, \qquad 1\le k\le \frac{T}{\tau}-1.
\ee
In view of the assumption (A), noting \eqref{p_k} and \eqref{nonlbd}, we get
\be\label{r_bound}
|r_j^k|\lesssim |\hat{e}_j^{\ep,k+1}|+|\hat{e}_j^{\ep,k}|+|\hat{e}_j^{\ep,k-1}|+
|\hat{f}^{\ep,k+1}_j|+|\hat{f}^{\ep,k-1}_j|,\qquad j\in \mathcal{T}_M. \ee
This implies that
\be\label{es1_r}
|(r^k,\Lbrack\hat{e}^\ep\Rbrack^k)|\lesssim  \mathcal{A}^k+\mathcal{A}^{k-1},\qquad
1\le k\le \frac{T}{\tau}-1.
\ee
By Cauchy inequality, we have
\be\label{es1_l}
\left|\mathrm{Im}\big(\hat{\xi}^{\ep,k},
\Lbrack\hat{e}^\ep\Rbrack^k\big)\right|\lesssim
\|\hat{\xi}^{\ep,k}\|^2+\|\hat{e}^{\ep,k+1}\|^2+\|\hat{e}^{\ep,k-1}\|^2
\lesssim \|\hat{\xi}^{\ep,k}\|^2+\mathcal{A}^k+\mathcal{A}^{k-1}.
\ee
In view of \eqref{n_eq}, \eqref{r_bound} and \eqref{Gp}, and using Cauchy inequality, we find
\bea\label{es2_q}
&&\left|\mathrm{Re}(H^{\ep,k}\Lbrack\hat{e}^\ep\Rbrack^k+\hat{\xi}^{\ep,k},
\tau\dt_t\hat{e}^{\ep,k})\right|\nn\\
&&=\tau\,\left|\mathrm{Im}\big(H^{\ep,k}\Lbrack\hat{e}^\ep\Rbrack^k+\hat{\xi}^{\ep,k},
(-\dt_x^2+H^{\ep,k})\Lbrack\hat{e}^\ep\Rbrack^k+r^k+\hat{\xi}^{\ep,k}\big)\right|\nn\\
&&\lesssim\tau\left(1+\|H^{\ep,k}\|_\infty+\|\dt_x^+H^{\ep,k}\|_\infty\right)\left(
\|\hat{\xi}^{\ep,k}\|^2+\|\dt_x^+\hat{\xi}^{\ep,k}\|^2+\mathcal{A}^k+\mathcal{A}^{k-1}\right)\nn\\
&&\lesssim\tau\left(\|\hat{\xi}^{\ep,k}\|^2+\|\dt_x^+\hat{\xi}^{\ep,k}\|^2
+\mathcal{A}^k+\mathcal{A}^{k-1}\right), \qquad
1\le k\le \frac{T}{\tau}-1.
\eea
Combining  \eqref{q1} and \eqref{nonlbd}, we have
\bea\label{q_1}
&&\left|\mathrm{Re}(q_1^{k},4\tau\dt_t\hat{e}^{\ep,k})\right|=4\tau\,
\left|\mathrm{Im}\left(q_1^k,(-\dt_x^2+H^{\ep,k})
\Lbrack\hat{e}^{\ep}\Rbrack^k+r^k+\hat{\xi}^{\ep,k}\right)\right| \nn\\
&&\lesssim \tau(1+\|H^{\ep,k}\|_\infty)(\|\dt_x^+\hat{e}^{\ep,k+1}\|^2+\|\dt_x^+\hat{e}^{\ep,k-1}\|^2
+\|\dt_x^+q_1^k\|^2+\|q_1^k\|^2\nn\\
&&\quad +\|r^k\|^2+\|\hat{\xi}^{\ep,k}\|^2
+\|\hat{e}^{\ep,k+1}\|^2+\|\hat{e}^{\ep,k-1}\|^2)\nn\\
&&\lesssim \tau(\|\hat{\xi}^{\ep,k}\|^2 +\mathcal{A}^k+\mathcal{A}^{k-1}), \qquad
1\le k\le \frac{T}{\tau}-1.
\eea
In view of \eqref{q1}, we get
\bea
\mathrm{Re}(q_2^k,4\tau\dt_t\hat{e}^{\ep,k})
&=&2\,\mathrm{Re}\big(g(\hat{E}^{\ep,k+1},\hat{E}^{\ep,k-1})
(\Lbrack\hat{f}^\ep\Rbrack^k-p^k),E^\ep(\cdot,t_{k+1})-E^\ep(\cdot,t_{k-1})\big)\nonumber\\
&&-\big(\Lbrack\hat{f}^\ep\Rbrack^k-p^k,2\tau\dt_t (\gamma_{_B}(|\hat{E}^{\ep,k}|^2))\big)
=q^k+\big(\Lbrack\hat{f}^\ep\Rbrack^k-p^k,2\tau\dt_t p^k\big),
\eea
where
\begin{equation*}
q^k=2\,\mathrm{Re}\big((g(\hat{E}^{\ep,k+1},\hat{E}^{\ep,k-1})-
\Lparen E^\ep\Rparen(\cdot,t_{k}))(\Lbrack\hat{f}^\ep\Rbrack^k-p^k),
E^\ep(\cdot,t_{k+1})-E^\ep(\cdot,t_{k-1})\big).
\end{equation*}
By Assumption (A) and \eqref{nonlbd}, we have
$$
|q^k|\lesssim \tau\|\p_t E^\ep\|_{L^\infty}(\mathcal{A}^k+\mathcal{A}^{k-1})\lesssim \tau(\mathcal{A}^k+\mathcal{A}^{k-1}), \qquad 1\le k\le \frac{T}{\tau}-1.
$$
Combining the above inequalities, we obtain
\be\label{es2_r}
\left|4\,\mathrm{Re}\left(r^k,\tau\dt_t\hat{e}^{\ep,k}\right)
-\left(\Lbrack\hat{f}^\ep\Rbrack^k-p^k,2\tau\dt_t p^k\right)\right|\lesssim\tau(\|\hat{\xi}^{\ep,k}\|^2+\mathcal{A}^k+\mathcal{A}^{k-1}).
\ee
Hence it can be concluded from \eqref{sp}, \eqref{es1_r}, \eqref{es1_l},
\eqref{es2_q} and \eqref{es2_r} that 
\bea\label{sp2}
&&\mathcal{A}^k-\mathcal{A}^{k-1}-(p^k,p^{k+1}-p^{k-1})-\tau(\hat{\eta}^{\ep,k},\hat{u}^{\ep,k+1/2}+
\hat{u}^{\ep,k-1/2})\nonumber\\
&&\lesssim \tau(\|\hat{\xi}^{\ep,k}\|^2+\|\dt_x^+\hat{\xi}^{\ep,k}\|^2
+\mathcal{A}^k+\mathcal{A}^{k-1}), \qquad 1\le k\le \frac{T}{\tau}-1.
\eea
Summing the above equation for $k=1,2,\cdots,m\le\fl{T}{\tau}-1$ and noting
$p^0={\bf 0}$ in (\ref{pjk234}),
we have
\bea\label{sp3}
&&\mathcal{A}^m-\mathcal{A}^0 -(p^m,p^{m+1})-\tau\sum\limits_{l=1}^m(\hat{\eta}^{\ep,l},\hat{u}^{\ep,l+1/2}+
\hat{u}^{\ep,l-1/2})\nonumber\\
&&\lesssim \tau \mathcal{A}^0+ \tau\sum\limits_{l=1}^m
(\|\hat{\xi}^{\ep,l}\|^2+\|\dt_x^+\hat{\xi}^{\ep,l}\|^2+
\mathcal{A}^l), \qquad 1\le m\le \frac{T}{\tau}-1.
\eea
Noting \eqref{innpX_M} and using Sobolev and Cauchy inequalities, we obtain
\bea\label{u_e}
&&-\fl{\mathcal{A}^m}{4}+\tau\sum\limits_{l=1}^m\left(\hat{\eta}^{\ep,l},\hat{u}^{\ep,l+1/2}+
\hat{u}^{\ep,l-1/2}\right)\nonumber\\
&&=-\fl{\mathcal{A}^m}{4}+\sum\limits_{l=1}^m\left((-\dt_x^2)^{-1}\hat{\eta}^{\ep,l},
\hat{f}^{\ep,l+1}-p^{l+1}-(\hat{f}^{\ep,l-1}-p^{l-1})\right)\nn\\
&&=-\fl{\mathcal{A}^m}{4}-2\tau\sum\limits_{l=2}^{m-1}\left(\dt_t(-\dt_x^2)^{-1}\hat{\eta}^{\ep,l},
\hat{f}^{\ep,l}-p^{l}\right)\nn\\
&&\quad +\sum\limits_{l=m}^{m+1}\left((-\dt_x^2)^{-1}\hat{\eta}^{\ep,l-1},\hat{f}^{\ep,l}-p^{l}\right)
-\sum\limits_{l=0}^1\left((-\dt_x^2)^{-1}\hat{\eta}^{\ep,l+1},\hat{f}^{\ep,l}-p^{l}\right)\nn\\
&&\lesssim \mathcal{A}^0+
\tau\sum\limits_{l=2}^{m-1} (\|\dt_t\hat{\eta}^{\ep,l}\|^2+\mathcal{A}^l)+
\sum\limits_{l=1}^2\|\hat{\eta}^{\ep,l}\|^2
+\sum\limits_{l=m-1}^{m}\|\hat{\eta}^{\ep,l}\|^2.
\eea
Recalling that
\be\label{pmpmp1}
\left(p^m,p^{m+1}\right)\le \fl{C_{_B}}{2}(\|\hat{e}^{\ep,m}\|^2+
\|\hat{e}^{\ep,m+1}\|^2)\le \fl{1}{2}\mathcal{A}^m,\qquad 1\le m\le \frac{T}{\tau}-1.
\ee
Combining (\ref{sp3}), \eqref{u_e} and \eqref{pmpmp1}, there exists $0<\tau_1\le \frac{1}{16}$ such that when $0<\tau\le \tau_1$, we have
\bea  \label{Amin}
\mathcal{A}^m&\lesssim&\mathcal{A}^0+ \tau\sum\limits_{l=1}^{m-1}\mathcal{A}^l
+\sum\limits_{l=1}^2\|\hat{\eta}^{\ep,l}\|^2
+\sum\limits_{l=m-1}^{m}\|\hat{\eta}^{\ep,l}\|^2\nn\\
&&+\tau\sum\limits_{l=1}^m
(\|\hat{\xi}^{\ep,l}\|^2+\|\dt_x^+\hat{\xi}^{\ep,l}\|^2)+
\tau\sum\limits_{l=2}^{m-1}\|\dt_t\hat{\eta}^{\ep,l}\|^2,
\qquad 1\le m\le \frac{T}{\tau}-1.
\eea
By Lemma \ref{initial-e} and using the discrete Sobolev inequality, we have
\be\label{u_e1}
\ep\|\dt_x^+\hat{u}^{\ep,1/2}\|\lesssim \ep\|\dt_t^+(\hat{f}^{\ep,0}- p^0)\|\lesssim \ep\|\dt_t^+\hat{f}^{\ep,0}\|+\ep\|\dt_t^+ \hat{e}^{\ep,0}\|
\lesssim \fl{\tau^2}{\ep},
\ee
which together with Lemma \ref{initial-e} yields that
\be\label{A0bd}
\mathcal{A}^0\lesssim \left(h^2+\fl{\tau^2}{\ep}\right)^2.
\ee
Plugging (\ref{A0bd}) into (\ref{Amin}) and noting Lemma \ref{local_e}, we get
\be
\mathcal{A}^m\lesssim \left(h^2+\fl{\tau^2}{\ep}\right)^2+\tau\sum\limits_{l=1}^{m-1}\mathcal{A}^l,\qquad
1\le m\le \frac{T}{\tau}-1.
\ee
Applying the discrete Gronwall inequality, when $0<\tau\le \tau_1$, we obtain
$$\mathcal{A}^m\lesssim \left(h^2+\fl{\tau^2}{\ep}\right)^2, \qquad
0\le m\le \frac{T}{\tau}-1,$$
which completes the proof of Theorem \ref{thm2} by noting (\ref{energy_d}).
\hfill $\square$

\bigskip

\begin{theorem}\label{thm3}
Under the assumptions (A)-(C), there exists $\tau_2>0$
sufficiently small and independent of $0<\ep\le1$, when $0<\tau\le \tau_2$ and $0<h\le \frac{1}{2}$,
we have the following error estimate of the
scheme \eqref{scheme}
\be\label{errort2}
\|\hat{e}^{\ep,k}\|+\|\dt_x^+\hat{e}^{\ep,k}\|+\|\hat{f}^{\ep,k}\|\lesssim h^2+\tau^2+\tau\ep^{\alpha^*}+\ep^{1+\alpha^*}, \quad 0\le k\le \fl{T}{\tau}.
\ee
\end{theorem}

 Define another set of error functions $\wt{e}^{\ep,k}\in X_M$ and $\wt{f}^{\ep,k}\in X_M$ as
\be
\wt{e}^{\ep,k}_j=\wt{E^\ep}(x_j,t_k)-\hat{E}_j^{\ep,k},\quad
\wt{f}_j^{\ep,k}=-\hat{F}_j^{\ep,k},\quad j\in \mathcal{T}_M^0,\quad 0\le k\le \fl{T}{\tau},
\ee
where  $\wt{E^\ep}$ is the solution of the NLSE-OP \eqref{ANLS1D}, and their corresponding
local truncation errors $\wt{\xi}^{\ep,k}\in X_M$ and $\wt{\eta}^{\ep,k}\in X_M$  as
\be\label{local2}
\begin{split}
\wt{\xi}_j^{\ep,k}&=i\dt_t \wt{E}^\ep(x_j,t_k)+(\dt_x^2-H_j^{\ep,k})\Lparen\wt{E}^\ep\Rparen(x_j,t_k)\\
&\quad+\gamma_{_B}(|\wt{E}^\ep(x_j,t_k)|^2)
g\big(\wt{E}^\ep(x_j,t_{k+1}),\wt{E}^\ep(x_j,t_{k-1})\big)\\
&=i\dt_t \wt{E}^\ep(x_j,t_k)+\big(\dt_x^2+|\wt{E}^\ep(x_j,t_k)|^2-H_j^{\ep,k}\big)
\Lparen\wt{E}^\ep\Rparen(x_j,t_k),\\
\wt{\eta}_j^{\ep,k}&=-\ep^2\dt_t^2\gamma_{_B}(|\wt{E}^\ep(x_j,t_k)|^2)
=-\ep^2\dt_t^2(|\wt{E}^\ep(x_j,t_k)|^2), \qquad j\in \mathcal{T}_M.
\end{split}
\ee
\begin{lemma}\label{local2e}
Under the assumption (A), when $0<h\le \frac{1}{2}$ and $0<\tau\le \frac{1}{2}$, we have
\be
\|\wt{\xi}^{\ep,k}\|+\|\dt_x^+\wt{\xi}^{\ep,k}\|\lesssim h^2+\tau^2+\tau\ep^{\alpha^*},\quad
\|\wt{\eta}^{\ep,k}\|\lesssim \ep^2,\quad \|\dt_t\wt{\eta}^{\ep,k}\|\lesssim \ep^{1+\alpha^*}.
\ee
\end{lemma}
\emph{Proof.} Similar to the proof of Lemma \ref{local_e}, we can get that
\begin{align*} \wt{\xi}^{\ep,k}_j=&\fl{h^2}{12}\int_0^1(1-s)^3\sum\limits_{m=\pm1}\sum\limits_{l=\pm1}
\wt{E}^\ep_{xxxx}(x_j+slh,t_k+m \tau)ds\\
&-\fl{\tau^2}{4}\int_{0}^1(1-s)^2
\sum\limits_{m=\pm 1}\p_{tt}(\wt{E}_{xx}^\ep+|\wt{E}^\ep|^2\wt{E}^\ep)(x_j,t_k+ms \tau)ds\\
&+\fl{\tau^2}{2}\int_0^1(1-s)\sum\limits_{m=\pm1}
\wt{E}^\ep_{xxtt}(x_j,t_k+ms\tau)ds+A_2\\
&+\fl{\tau^2}{2}\big(|\wt{E}^\ep(x_j,t_{k})|^2-H_j^{\ep,k}\big)
\int_0^1(1-s)\sum\limits_{m=\pm1}\wt{E}^\ep_{tt}(x_j,t_k+ms\tau)ds,
\end{align*}
where
\begin{align*}
|A_2|&=\left|\fl{1}{2\tau}\int_{-\tau}^\tau \wt{E}^\ep(x_j,t_k+s) G^\ep\left(x_j,\fl{t_k+s}{\ep}\right)ds-\wt{E}^\ep(x_j,t_k)H_j^{\ep,k}\right|\\
&=\left|\fl{\tau}{2}\int_{-1}^1G^\ep\left(x_j,\fl{t_k+s\tau}{\ep}\right)
\int_0^s\wt{E}^\ep_t(x_j,t_k+\tht\tau)d\tht ds\right|\\
&\lesssim \tau \|G^\ep\|_{L^\infty}\,\|\wt{E}^\ep_t\|_{L^\infty}
\lesssim \tau\ep^{\alpha^*}, \qquad j\in \mathcal{T}_M, \quad 1\le k\le \frac{T}{\tau}-1.
\end{align*}
Recalling \eqref{ANLS}, \eqref{Gp} and assumption (A), and using integration by parts,  we have
\bea \label{tp}
&&\tau^2\left|\int_0^1(1-s)\wt{E}^\ep_{tt}(x_j,t_k+s\tau)ds\right|\nn\\
&&=\left|\tau^2\int_0^1(1-s)\big(\wt{E}^\ep_{xxt}+(|\wt{E}^{\ep}|^2\wt{E}^\ep
)_t\big)(x_j,t_k+s\tau)ds\right.\nn\\
&&\quad \left.-\tau^2\int_0^1(1-s)\p_s\left[\wt{E}^\ep(x_j,s) G^\ep\left(x_j,\fl{s}{\ep}\right)\right]\Big|_{(t_k+s\tau)}ds\right|\nn\\
&&\le\tau \left|\wt{E}^\ep(x_j,t_k)G^\ep\left(x_j,\fl{t_k}{\ep}\right)-\int_0^1
\wt{E}^\ep(x_j,t_k+s\tau)G^\ep\left(x_j,\fl{t_k+s\tau}{\ep}\right)ds\right|
\nn\\
&&\quad+\tau^2\left|\int_0^1(1-s)\big(\wt{E}^\ep_{xxt}+(|\wt{E}^{\ep}|^2\wt{E}^\ep)_t
\big)(x_j,t_k+s\tau)ds\right|\nn\\
&&\lesssim  \tau^2+\tau\ep^{\alpha^*}, \qquad j\in \mathcal{T}_M, \quad 1\le k\le \frac{T}{\tau}-1.
\eea
Similarly, we can get that
\begin{align*}
&\tau^2\left|\int_{0}^1(1-s)^2\sum\limits_{m=\pm 1}\p_{tt}(\wt{E}_{xx}^\ep+|\wt{E}^\ep|^2\wt{E}^\ep)(x_j,t_k+ms \tau)ds\right|\lesssim \tau^2+\tau\ep^{\alpha^*},\\
&\tau^2\left|\int_0^1(1-s)\sum\limits_{m=\pm1}
\wt{E}^\ep_{xxtt}(x_j,t_k+ms\tau)ds\right|\lesssim \tau^2+\tau\ep^{\alpha^*}.
\end{align*}
Hence we can conclude that
$$\|\wt{\xi}^{\ep,k}\|\lesssim h^2+\tau^2+\tau\ep^{\alpha^*},\qquad
 1\le k\le \frac{T}{\tau}-1.$$
Similarly, we can get
\[ \|\dt_x^+\wt{\xi}^{\ep,k}\|\lesssim h^2+\tau^2+\tau\ep^{\alpha^*},\qquad
 1\le k\le \frac{T}{\tau}-1.\]
By assumption (A), it is easy to get that
\[
\begin{split}
&\left|\p_{tt}|\wt{E}^\ep(x,t)|^2\right|=\left|-2\,\mathrm{Im}\left(\overline{\wt{E}^\ep_t}\wt{E}^\ep_{xx}
+\overline{\wt{E}^\ep}\wt{E}^\ep_{xxt}\right)\right|\lesssim 1,\qquad x\in\Omega, \quad 0\le t\le T,\\
&\left|\p_{ttt}|\wt{E}^\ep(x,t)|^2\right|=\left| -2\,\mathrm{Im}\left(\overline{\wt{E}^\ep_{tt}}\wt{E}^\ep_{xx}
+2\overline{\wt{E}^\ep_t}\wt{E}^\ep_{xxt}+\overline{\wt{E}^\ep}
\wt{E}^\ep_{xxtt}\right)\right|\lesssim \ep^{\alpha^*-1},
\end{split}
\]
which indicate that
$$\|\wt{\eta}^{\ep,k}\|\lesssim \ep^2,\quad 1\le k\le \frac{T}{\tau}-1;\quad
\|\dt_t \wt{\eta}^{\ep,k}\|\lesssim \ep^{1+\alpha^*},\quad 2\le k\le \frac{T}{\tau}-2.$$
Thus the proof is completed.
\hfill $\square$

Analogous to Lemma \ref{initial-e}, we have error bounds of
$\wt{e}^{\ep,k}$, $\wt{f}^{\ep,k}$ at the first step.
\begin{lemma} \label{initial-e2}
Under the assumptions (A) and (B), when $0<h\le \frac{1}{2}$ and $0<\tau\le \frac{1}{2}$, we have
\begin{align*}
&\wt{e}_j^{\ep,0}=\wt{f}_j^{\ep,0}=0,\qquad  |\wt{e}_j^{\ep,1}|+|\dt_x^+\wt{e}_j^{\ep,1}|\lesssim \tau^2+\tau\ep^{\alpha^*},\\
&|\dt_t^+\wt{e}_j^{\ep,0}|\lesssim \tau+\ep^{\alpha^*}, \qquad |\wt{f}_j^{\ep,1}|\lesssim \tau^2,\qquad
|\dt_t^+\wt{f}_j^{\ep,0}|\lesssim \tau,\qquad j\in \mathcal{T}_M.
\end{align*}
\end{lemma}
\emph{Proof.} It follows from \eqref{Zak1d} and \eqref{ANLS} that
$\p_tE^\ep(x_j,0)=\p_t\wt{E}^\ep(x_j,0)=\phi_2(x_j)$ for $j\in \mathcal{T}_M^0$.
By \eqref{F_1}, \eqref{tp} and assumption (B), we get
$$
|\wt{e}_j^{\ep,1}|=\left|\tau^2\int_0^1(1-s)\wt{E}^\ep_{tt}(x_j,s\tau)ds-\fl{\tau^2}{2}
E^\ep_{tt}(x_j,0)\right|\lesssim \tau^2+\tau\ep^{\alpha^*}, \qquad j\in \mathcal{T}_M.$$
Similarly, we have
 \[ |\dt_x^+\wt{e}_j^{\ep,1}|\lesssim \tau^2+\tau\ep^{\alpha^*}, \qquad j\in \mathcal{T}_M.\]
Moreover, it is easy to get that
$$|\wt{f}^{\ep,1}_j|=|F_j^{\ep,1}|\lesssim \tau^2 |F^\ep_{tt}(x_j,0)|\lesssim \tau^2,
\qquad j\in \mathcal{T}_M.$$
The rest can be obtained similarly and details are omitted here for brevity.
\hfill $\square$

\bigskip

\emph{Proof of Theorem \ref{thm3}.}
Subtracting \eqref{scheme} from \eqref{local2}, we obtain the error equations
\be\label{2e_eq}
\begin{split}
&i\dt_t\wt{e}^{\ep,k}_j=(-\dt_x^2+H_j^{\ep,k})\Lbrack\wt{e}^\ep\Rbrack_j^k
+\wt{r}_j^k+\wt{\xi}_j^{\ep,k},\\
&\ep^2\dt_t^2\wt{f}_j^{\ep,k}=\dt_x^2\Lbrack\wt{f}^\ep\Rbrack_j^k
+\ep^2\dt_t^2\wt{p}_j^k+\wt{\eta}_j^{\ep,k},\qquad j\in \mathcal{T}_M, \quad 1\le k\le \frac{T}{\tau}-1,
\end{split}
\ee
where $\wt{r}^k\in X_M$ and $\wt{p}^k\in X_M$ defined as
\begin{align*}
\wt{r}_j^k&=-|\wt{E}^\ep(x_j,t_k)|^2\Lparen\wt{E}^\ep\Rparen(x_j,t_k)
+\left(\gamma_{_B}(|\hat{E}_j^{\ep,k}|^2)-\Lbrack\hat{F}^\ep\Rbrack^k_j\right) g(\hat{E}_j^{\ep,k+1},\hat{E}_j^{\ep,k-1}),\\
\wt{p}_j^k&=|\wt{E}^\ep(x_j,t_k)|^2-\gamma_{_B}(|\hat{E}_j^{\ep,k}|^2),\qquad j\in \mathcal{T}_M,
\quad 1\le k\le \frac{T}{\tau}-1.
\end{align*}
Let $\wt{u}^{\ep,k+\fl{1}{2}}\in X_M$ be the solution of the equation
$$-\dt_x^2\wt{u}_j^{\ep,k+\fl{1}{2}}=\dt_t^+(\wt{f}_j^{\ep,k}-\wt{p}_j^k),
\qquad j\in\mathcal{T}_M,\quad 0\le k\le \frac{T}{\tau}-1.$$
Define another discrete `energy'
\begin{align}
&\wt{\mathcal {A}}^k=C_{_B}(\|\wt{e}^{\ep,k}\|^2+\wt{e}^{\ep,k+1}\|^2)+
\|\dt_x^+\wt{e}^{\ep,k}\|^2+\|\dt_x^+\wt{e}^{\ep,k+1}\|^2\nonumber\\
&\qquad+\ep^2\|\dt_x^+\wt{u}^{\ep,k+1/2}\|^2
+\fl{1}{2}(\|\wt{f}^{\ep,k}\|^2+\|\wt{f}^{\ep,k+1}\|^2),
\qquad 0\le k\le \frac{T}{\tau}-1.\label{energy_d2}
\end{align}
Applying the same approach as in the proof of Theorem 3.1 and noting
$(\wt{p}^k,\wt{p}^{k+1})\le \fl{1}{2}\wt{\mathcal{A}}^k$, there exists
$0<\tau_2\le\frac{1}{16}$ sufficiently small and independent of $0<\ep\le1$ such that when $0<\tau\le\tau_2$,
\[ 
\wt{\mathcal{A}}^k\lesssim \wt{\mathcal{A}}^0+
\tau\sum\limits_{l=1}^{k-1}\wt{\mathcal{A}}^l
+\sum\limits_{l=1}^2\|\wt{\eta}^{\ep,l}\|^2
+\sum\limits_{l=k-1}^{k}\|\wt{\eta}^{\ep,l}\|^2\\
+\tau\sum\limits_{l=1}^k
(\|\wt{\xi}^{\ep,l}\|^2+\|\dt_x^+\wt{\xi}^{\ep,l}\|^2)+
\tau\sum\limits_{l=2}^{k-1}\|\dt_t\wt{\eta}^{\ep,l}\|^2.
\]
By Lemma \ref{initial-e2} and the discrete Sobolev inequality, we deduce that
$$\ep\|\dt_x^+\wt{u}^{\ep,1/2}\|\lesssim \ep\|\dt_t^+\wt{f}^{\ep,0}\|+\ep
\|\dt_t^+\wt{e}^{\ep,0}\|\lesssim \ep\tau+\ep^{1+\alpha^*},$$
which together with Lemma \ref{initial-e2} yields that
$$\wt{\mathcal{A}}^0\lesssim (\tau^2+\tau\ep^{\alpha^*}+\ep^{1+\alpha^*})^2.$$
By Lemma \ref{local2e}, when $0<\tau\le \tau_2$ and $0<h\le \frac{1}{2}$, we have
$$\wt{\mathcal {A}}^k\lesssim \left(h^2+\tau^2+\tau\ep^{\alpha^*}
+\ep^{1+\alpha^*}\right)^2+\tau\sum\limits_{l=1}^{k-1}\wt{\mathcal{A}}^l, \qquad
1\le k\le \frac{T}{\tau}-1.$$
Using the discrete Gronwall inequality, when $0<\tau\le \tau_2$, one has
$$\wt{\mathcal{A}}^k\lesssim \left(h^2+\tau^2+\tau\ep^{\alpha^*}+
\ep^{1+\alpha^*}\right)^2,\qquad 1\le k\le \frac{T}{\tau}-1.$$
Noting \eqref{energy_d2}, we get
$$\|\wt{e}^{\ep,k}\|+\|\dt_x^+\wt{e}^{\ep,k}\|+\|\wt{f}^{\ep,k}\|\lesssim
h^2+\tau^2+\tau\ep^{\alpha^*}+\ep^{1+\alpha^*}, \qquad 0\le k\le \frac{T}{\tau}.$$
Combining the above inequality and (\ref{Gp}), using the triangle inequality and
noting assumption (C),  we obtain
\[
\begin{split}
&\|\hat{e}^{\ep,k}\|+\|\dt_x^+\hat{e}^{\ep,k}\|\lesssim\|\wt{e}^{\ep,k}\|
+\|\dt_x^+\wt{e}^{\ep,k}\|+\|E^\ep(\cdot,t_k)-\wt{E}^\ep(\cdot,t_k)\|_{H^1}\\
&\qquad \qquad \qquad \quad\  \lesssim h^2+\tau^2+\tau\ep^{\alpha^*}+\ep^{1+\alpha^*},
\qquad 0\le k\le \frac{T}{\tau},\\
&\|\hat{f}^{\ep,k}\|\lesssim \|\wt{f}^{\ep,k}\|+\|F^\ep(\cdot,t_k)-F(\cdot,t_k)\|_{L^2}
\lesssim h^2+\tau^2+\tau\ep^{\alpha^*}+\ep^{1+\alpha^*},
\end{split}
\]
which complete the proof of Theorem \ref{thm3}.
\hfill $\square$

\bigskip

\emph{Proof of Theorem \ref{thm1}.}
When $0<\tau\le \min\left\{\frac{1}{16},\tau_1,\tau_2\right\}$ and $0<h\le \frac{1}{2}$,
combining (\ref{errort1}) and (\ref{errort2}), we have for $0\le k\le \frac{T}{\tau}$
\be
\|\hat{e}^{\ep,k}\|+\|\dt_x^+\hat{e}^{\ep,k}\|+\|\hat{f}^{\ep,k}\|
\lesssim h^2+\min_{0<\ep\le1}\left\{\tau^2+\ep^{\alpha^*}(\tau+\ep),
\frac{\tau^2}{\ep}\right\}\lesssim h^2+\tau^{1+\fl{\alpha^*}{2+\alpha^*}}.
\ee
This, together with the inverse inequality \cite{Thomee}, implies
\[
\|\hat{E}^{\ep,k}\|_\infty-\|E^\ep(\cdot,t_k)\|_\infty\le\|\hat{e}^{\ep,k}\|_\infty\lesssim \|\delta_x^+\hat{e}^{\ep,k}\|
\lesssim h^2+\tau^{1+\fl{\alpha^*}{2+\alpha^*}}, \qquad 0\le k\le \frac{T}{\tau}.
\]
Thus, there exist $h_1>0$ and $\tau_3>0$ sufficiently small and independent of $0<\ep\le1$
such that when $0<h\le h_1$ and $0<\tau\le \tau_3$,
\[ \|\hat{E}^{\ep,k}\|_\infty\le 1+\|E^\ep(\cdot,t_k)\|_\infty\le 1+M_0, \qquad 0\le k\le \frac{T}{\tau}.
\]
Taking $h_0=\min\left\{\frac{1}{2}, h_1\right\}$ and $\tau_0=\min\left\{\frac{1}{16}, \tau_1,\tau_2,\tau_3\right\}$,
when $0<h\le h_0$ and $0<\tau\le \tau_0$, the numerical method \eqref{scheme}
collapses to \eqref{scheme1}, i.e.
\[ \hat{E}^{\ep,k}_j=E^{\ep,k}_j,\qquad \hat{F}^{\ep,k}_j=F^{\ep,k}_j, \qquad j\in \mathcal{T}_M^0,
\quad 0\le k\le \frac{T}{\tau}.\]
Thus the proof is completed. \hfill $\square$

\begin{remark}
The error bounds in Theorem \ref{thm1} are still valid in high dimensions, e.g.
$d=2,3$, provided that an additional condition on the time step $\tau$ is added
\[\tau=o\left(C_d(h)^{1-\fl{\alpha^*}{2+2\alpha^*}}\right),\]
with
$$
 C_d(h)\sim
 \left\{
\begin{aligned}
&\fl{1}{|\mathrm{ln}\, h|},\quad &d=2,\\
&h^{1/2},\quad &d=3. \\
\end{aligned}
\right.
$$
The reason is due to the discrete Sobolev inequality \cite{Cai2012,Cai2013,Bao2013b}
$$\|\psi_h\|_\infty\le \fl{1}{C_d(h)}\|\psi_h\|_{H^1},
$$
where $\psi_h$ is a mesh function over $\Omega$ with homogeneous Dirichlet
boundary condition.
\end{remark}

\section{Numerical results}

In this section, we present numerical results for the
 ZS \eqref{Zak} by our proposed finite difference method.
In order to do so, we take $d=1$ in \eqref{Zak}
and the initial condition is taken as
$$E_0(x)=e^{-x^2/2},\quad \og_0(x)=e^{-x^2/4},\quad \og_1(x)=e^{-x^2/3}\sin(x), \qquad x\in {\mathbb R}.$$
We mainly consider two types of initial data 

Case I. well-prepared initial data, i.e., $\alpha=1$ and $\beta=0$;

Case II. ill-prepared initial data, i.e., $\alpha=0$ and $\beta=0$.

\medskip

In practical computation, the problem is truncated on a bounded interval $\Omega=[-200,200]$,
which is large enough such that the homogeneous Dirichlet boundary condition does not
introduce significant errors. In addition, we introduce the following error functions
$$e^\ep(t_k):=\|e^{\ep,k}\|+\|\dt_x^+e^{\ep,k}\|,\quad n^\ep(t_k):=\|N^\ep(\cdot,t_k)-N^{\ep,k}\|,
\qquad k\ge0,$$
where $e^{\ep,k}_j=E^\ep(x_j,t_k)-E^{\ep,k}_j$ and $N^{\ep,k}_j=-|E^{\ep,k}_j|^2+F^{\ep,k}_j+
G^\ep(x_j,t_k/\ep)$ for $0\le j\le M$.
The ``exact" solution is obtained by the time splitting spectral method \cite{Bao2005} with very small mesh size $h=1/64$ and time step $\tau=10^{-6}$.

Table \ref{ill-h} depicts the spatial errors at $t=1$ with a fixed time step
$\tau=10^{-5}$ and Case II initial data
for different mesh size $h$ and $0<\ep\le1$.
It clearly demonstrates that our new finite difference method is
uniformly second order accurate in space for all $\ep \in (0,1]$.
The results for other initial data are analogous, e.g. different
$\alpha\ge0$ and $\beta\ge0$ and thus are omitted for brevity.

\begin{table}[t!]
\def\temptablewidth{1\textwidth}
\vspace{-12pt}
\caption{Spatial error analysis at time $t=1$ for Case II, i.e.
$\alpha=\beta=0$.}\label{ill-h}
{\rule{\temptablewidth}{1pt}}
\begin{tabular*}{\temptablewidth}{@{\extracolsep{\fill}}cccccccc}
$e^\ep(1)$&$h_0=0.2$&$h_0/2$&$h_0/2^2$&$h_0/2^3$&$h_0/2^4$&$h_0/2^5$\\\hline
$\ep=1$&	2.83E-2&	7.27E-3&	1.82E-3&	4.56E-4&	1.14E-4& 2.85E-5\\
rate&-&	1.96&	1.99&	2.00&	2.00& 2.00\\\hline
$\ep=1/2$&	2.62E-2&	6.73E-3&	1.69E-3&	4.23E-4&	1.06E-4& 2.65E-5\\
rate&-&	1.96&	1.99&	2.00&	2.00& 2.00\\\hline
$\ep=1/2^{2}$&	2.52E-2&	6.44E-3&	1.61E-3&	4.03E-4&	1.01E-4& 2.53E-5\\
rate&-&	1.97&	2.00&	2.00&	2.00& 2.00\\\hline
$\ep=1/2^{3}$&	2.63E-2&	6.73E-3&	1.69E-3&	4.23E-4&	1.06E-4& 2.65E-5\\
rate&-&	1.97&	1.99&	2.00&	2.00& 2.00\\\hline
$\ep=1/2^{4}$&	2.64E-2&	6.67E-3&	1.67E-3&	4.18E-4&	1.05E-4& 2.63E-5\\
rate&-&	1.98&	2.00&	2.00&	2.00& 2.00\\\hline
$\ep=1/2^{5}$&	2.68E-2&	6.80E-3&	1.70E-3&	4.26E-4&	1.07E-4& 2.68E-5\\
rate&-&	1.98&	2.00&	2.00&	2.00& 2.00\\\hline
$\ep=1/2^{6}$&	2.69E-2&	6.83E-3&	1.71E-3&	4.28E-4&	1.07E-4& 2.68E-5 \\
rate&-&	1.98&	2.00&	2.00&	2.00& 2.00\\
\bottomrule
  \bottomrule
$n^\ep(1)$&$h_0=0.2$&$h_0/2$&$h_0/2^2$&$h_0/2^3$&$h_0/2^4$&$h_0/2^5$\\\hline
$\ep=1$&	7.24E-3&	1.80E-3&	4.50E-4&	1.12E-4&	2.81E-5& 7.03E-6\\
rate&-&	2.01&	2.00&	2.00&	2.00& 2.00\\\hline
$\ep=1/2$&	9.52E-3&	2.36E-3&	5.90E-4&	1.47E-4&	3.69E-5& 9.23E-6\\
rate&-&	2.01&	2.00&	2.00&	2.00& 2.00\\\hline
$\ep=1/2^{2}$&	7.20E-3&	1.80E-3&	4.49E-4&	1.12E-4&	2.81E-5& 7.03E-6\\
rate&-&	2.00&	2.00&	2.00&	2.00& 2.00\\\hline
$\ep=1/2^{3}$&	4.76E-3&	1.18E-3&	2.95E-4&	7.36E-5&	1.84E-5& 4.60E-6\\
rate&-&	2.01&	2.00&	2.00&	2.00& 2.00\\\hline
$\ep=1/2^{4}$&	4.55E-3&	1.13E-3&	2.81E-4&	7.01E-5&	1.75E-5& 4.38E-6\\
rate&-&	2.02&	2.00&	2.00&	2.00& 2.00\\\hline
$\ep=1/2^{5}$&	4.52E-3&	1.12E-3&	2.78E-4&	6.95E-5&	1.74E-5& 4.35E-6\\
rate&-&	2.02&	2.00&	2.00&	2.00& 2.00\\\hline
$\ep=1/2^{6}$&	4.51E-3&	1.11E-3&	2.78E-4&	6.94E-5&	1.74E-5& 4.35E-6\\
rate&-&	2.02&	2.00&	2.00&	2.00& 2.00\\
\end{tabular*}
{\rule{\temptablewidth}{1pt}}
\end{table}

Table \ref{well-t}  presents the temporal errors at
$t=1$ with a fixed mesh size $h=2.5\times 10^{-4}$ and
Case I initial data for different time step $\tau$ and $0<\ep\le1$,
and respectively, Table \ref{ill-t} depicts similar results for Case II initial data.

\begin{table}[t!]
\def\temptablewidth{1\textwidth}
\setlength{\tabcolsep}{3pt}
\vspace{-12pt}
\caption{Temporal error analysis at time $t=1$ for Case I, i.e. $\alpha=1$ and $\beta=0$.}\label{well-t}
{\rule{\temptablewidth}{1pt}}
\begin{tabular*}{\temptablewidth}{@{\extracolsep{\fill}}ccccccccccc}
$e^\ep(1)$&$\tau_0=0.1$&$\tau_0/2$&$\tau_0/2^2$&$\tau_0/2^3$&$\tau_0/2^4$&$\tau_0/2^5$&
$\tau_0/2^6$&$\tau_0/2^7$\\\hline
$\ep=1$&	1.19E-1&	4.47E-2&	1.65E-2&	4.83E-3&	1.25E-3&	3.16E-4&   7.90E-5&   1.99E-5\\
rate&-&	1.42&	1.44&	1.77&	1.95&	1.99&  2.00& 1.99\\\hline
$\ep=1/2$&	7.80E-2&	3.66E-2&	1.46E-2&	4.33E-3&	1.12E-3&	2.83E-4&	7.10E-5&	1.79E-5\\
rate&-&	1.09&	1.33&	1.76&	1.94&	1.99&	2.00&	1.99\\\hline
$\ep=1/2^{2}$&	7.18E-2&	3.19E-2&	1.27E-2&	3.86E-3&	1.01E-3&	2.55E-4&	6.39E-5&	1.61E-5\\
rate&-&	1.17&	1.32&	1.72&	1.94&	1.99&	2.00&	1.99\\\hline
$\ep=1/2^{3}$&	7.12E-2&	3.67E-2&	1.35E-2&	3.80E-3&	9.79E-4&	2.47E-4&	6.19E-5&	1.56E-5\\
rate&-&	0.96&	1.45&	1.83&	1.96&	1.99&	2.00&	1.99\\\hline
$\ep=1/2^{4}$&	6.99E-2&	3.63E-2&	1.35E-2&	3.84E-3&	9.89E-4&	2.49E-4&	6.24E-5&	1.57E-5\\
rate&-&	0.94&	1.43&	1.81&	1.96&	1.99&	2.00&	1.99\\\hline
$\ep=1/2^{5}$&	6.97E-2&	3.63E-2&	1.36E-2&	3.87E-3&	9.96E-4&	2.51E-4&	6.29E-5&   1.58E-5\\
rate&-&	0.95&	1.42&	1.81&	1.96&	1.99&	2.00&	1.99\\\hline
$\ep=1/2^{6}$&6.97E-2	&3.63E-2&1.36E-2&3.87E-3&9.95E-4&2.50E-4&6.26E-5	&1.58E-5 \\
rate&-&	0.94&	1.42&	1.81&	1.96&	1.99&	2.00& 1.99\\
\toprule
\bottomrule
$n^\ep(1)$&$\tau_0=0.1$&$\tau_0/2$&$\tau_0/2^2$&$\tau_0/2^3$&$\tau_0/2^4$&$\tau_0/2^5$&
$\tau_0/2^6$&$\tau_0/2^7$\\\hline
$\ep=1$&	1.09E-2&	2.89E-3&	7.47E-4&	1.90E-4&	4.79E-5&	1.21E-5&     3.03E-6&   7.63E-7\\
rate&-&	1.91&	1.95&	1.98&	1.99&	1.99&	2.00&	1.99\\\hline
$\ep=1/2$&	2.15E-2&	6.09E-3&	1.59E-3&	4.05E-4&	1.02E-4&	2.56E-5&	6.46E-6&	1.63E-6\\
rate&-&	1.82&	1.94&	1.97&	1.99&	1.99&	1.99&	1.99\\\hline
$\ep=1/2^{2}$&	3.08E-2&	1.33E-2&	4.23E-3&	1.12E-3&	2.86E-4&	7.18E-5&	1.80E-5&	4.53E-6\\
rate&-&	1.22&	1.65&	1.91&	1.98&	1.99&	1.99&	1.99\\\hline
$\ep=1/2^{3}$&	1.62E-2&	6.76E-3&	2.72E-3&	1.23E-3&	3.57E-4&	9.16E-5&	2.30E-5&	5.78E-6\\
rate&-&	1.26&	1.31&	1.14&	1.79&	1.96&	1.99&	2.00\\\hline
$\ep=1/2^{4}$&	5.87E-3&	3.95E-3&	2.38E-3&	9.46E-4&	4.14E-4&	1.40E-4&	3.63E-5&	9.14E-6\\
rate&-&	0.57&	0.73&	1.33&	1.19&	1.57&	1.94&	1.99\\\hline
$\ep=1/2^{5}$&	7.39E-3&	1.98E-3&	1.08E-3&	7.62E-4&	4.01E-4&	1.45E-4&	6.20E-5&1.73E-5\\
rate&-&	1.90&	0.87&	0.51&	0.93&	1.47&	1.22&	1.84\\\hline
$\ep=1/2^{6}$&7.63E-3&2.72E-3&6.17E-4& 2.87E-4&2.24E-4& 1.49E-4&5.91E-5	&2.50E-5\\
rate&-&	1.49&	2.14&	1.10&	0.36&	0.59&	1.33& 1.24\\
\end{tabular*}
{\rule{\temptablewidth}{1pt}}
\end{table}

\begin{table}[t!]
\def\temptablewidth{1\textwidth}
\setlength{\tabcolsep}{3pt}
\vspace{-12pt}
\caption{Temporal error analysis at time $t=1$ for Case II, i.e. $\alpha=\beta=0$.}\label{ill-t}
{\rule{\temptablewidth}{1pt}}
\begin{tabular*}{\temptablewidth}{@{\extracolsep{\fill}}ccccccccccc}
$e^\ep(1)$&$\tau_0=0.1$&$\tau_0/2$&$\tau_0/2^2$&$\tau_0/2^3$&$\tau_0/2^4$&$\tau_0/2^5$&
$\tau_0/2^6$&$\tau_0/2^7$\\\hline
$\ep=1$&	1.19E-1&	4.47E-2&	1.65E-2&	4.83E-3&	1.25E-3&	3.16E-4&   7.90E-5&   1.99E-5\\
rate&-&	1.42&	1.44&	1.77&	1.95&	1.99&  2.00& 1.99\\\hline
$\ep=1/2$&	1.10E-1&	4.23E-2&	1.60E-2&	4.74E-3&	1.23E-3&	3.11E-4&   7.78E-5&   1.96E-5\\
rate&-&	1.38&	1.40&	1.76&	1.94&	1.99& 2.00& 1.99\\\hline
$\ep=1/2^{2}$&	1.03E-1&	4.13E-2&	1.50E-2&	4.52E-3&	1.18E-3&	2.98E-4&   7.45E-5&   1.88E-5\\
rate&-&	1.33&	1.46&	1.73&	1.94&	1.99&2.00 &1.99\\\hline
$\ep=1/2^{3}$&	7.87E-2&	4.08E-2&	1.57E-2&	4.65E-3&	1.21E-3&	3.07E-4&   7.68E-5&   1.93E-5\\
rate&-&	0.95&	1.38&	1.75&	1.94&	1.99&2.00 &1.99\\\hline
$\ep=1/2^{4}$&	7.05E-2&	3.66E-2&	1.35E-2&	3.81E-3&	9.77E-4&	2.46E-4&   6.15E-5&   1.55E-5\\
rate&-&	0.95&	1.44&	1.83&	1.96&	1.99&2.00&1.99 \\\hline
$\ep=1/2^{5}$&	7.05E-2&	3.59E-2&	1.35E-2&	3.87E-3&	9.98E-4&	2.51E-4&   6.28E-5&   1.58E-5\\
rate&-&	0.97&	1.41&	1.81&	1.96&	1.99&2.00  &1.99\\\hline
$\ep=1/2^{6}$&	6.99E-2&	3.62E-2&	1.35E-2&	3.86E-3&	9.96E-4&      2.51E-4&     6.29E-5&   1.58E-5\\
rate&-&	0.95&	1.42&	1.81&	1.95&1.99&  2.00&1.99\\
\toprule
\bottomrule
$n^\ep(1)$&$\tau_0=0.1$&$\tau_0/2$&$\tau_0/2^2$&$\tau_0/2^3$&$\tau_0/2^4$&$\tau_0/2^5$&
$\tau_0/2^6$&$\tau_0/2^7$\\\hline
$\ep=1$&	1.09E-2&	2.89E-3&	7.47E-4&	1.90E-4&	4.79E-5&	1.21E-5&     3.03E-6&   7.63E-7\\
rate&-&	1.91&	1.95&	1.98&	1.99&	1.99&	2.00&	1.99\\\hline
$\ep=1/2$&	2.62E-2&	7.37E-3&	1.93E-3&	4.90E-4&	1.24E-4&	3.11E-5&     7.78E-6&   1.96E-6\\
rate&-&	1.83&	1.94&	1.97&	1.99&	1.99&	2.00&	1.99\\\hline
$\ep=1/2^{2}$&	3.46E-2&	1.38E-2&	4.30E-3&	1.14E-3&	2.90E-4&	7.29E-5&     1.82E-5&   4.58E-6\\
rate&-&	1.33&	1.68&	1.91&	1.98&	1.99&	2.00&	1.99\\\hline
$\ep=1/2^{3}$&	2.21E-2&	1.07E-2&	3.71E-3&	1.31E-3&	3.70E-4&	9.45E-5&	2.38E-5&	5.96E-6\\
rate&-&	1.05&	1.52&	1.50&	1.83&	1.97&	1.99&	2.00\\\hline
$\ep=1/2^{4}$&	6.28E-3&	5.16E-3&	3.70E-3&	1.61E-3&	5.12E-4&	1.53E-4&	3.95E-5&	9.91E-6\\
rate&-&	0.28&	0.48&	1.20&	1.66&	1.74&	1.96&	1.99\\\hline
$\ep=1/2^{5}$&	5.73E-3&	2.34E-3&	1.34E-3&	1.11E-3&	6.74E-4&	2.25E-4&	7.10E-5& 1.91E-5\\
rate&-&	1.29&	0.80&	0.27&	0.72&	1.58&	1.67& 1.90\\\hline
$\ep=1/2^{6}$&	7.63E-3&	4.44E-3&	8.36E-4&	3.48E-4&	2.98E-4&	2.35E-4&	1.04E-4& 3.18E-5\\
rate&-&	0.78&	2.41&	1.26&	0.22&	0.35&	1.17& 1.71\\
\end{tabular*}
{\rule{\temptablewidth}{1pt}}
\end{table}

From Tables \ref{well-t} \& \ref{ill-t}, we can see that our numerical method
is `essentially' second-order in time for any fixed $0<\ep\le 1$
for both well-prepared and ill-prepared initial data. In fact, for each fixed $0<\ep\le 1$,
second order convergence in time is observed for $0<\tau\tau_0$ with $\tau_0>0$ independent
of $\ep$ except a small resonance region (cf. each row in Tables \ref{well-t} \& \ref{ill-t}), 
e.g. at $\tau=O(\ep^{3/2})$ for the well-prepared initial data Case I
and at $\tau=O(\ep)$ for the ill-prepared initial data Case II.
In fact, for well-prepared initial data Case I, in the resonance region
$\tau=O(\ep^{3/2})$, the convergence rate is downgraded to $4/3$;
and respectively, for the ill-prepared initial data Case II,
it is downgraded to first order in the resonance region $\tau=O(\ep)$;
which are listed in  Table \ref{de-w-t}. All these numerical results
demonstrate that our error bounds are sharp.

\begin{table}[h!]
\def\temptablewidth{1\textwidth}
\setlength{\tabcolsep}{3pt}
\vspace{-12pt}
\caption{Temporal error analysis at time $t=1$
for well-prepared and ill-prepared initial data in the resonance regions
with different $\tau$ and $\ep$.
}\label{de-w-t}
{\rule{\temptablewidth}{1pt}}
\begin{tabular*}{\temptablewidth}{@{\extracolsep{\fill}}ccccccccc}
Case I ($\tau=O(\ep^{3/2})$) &$\ep_0=1/2, \tau_0=0.1$
&$\ep_0/2^2, \tau_0/2^3$
&$\ep_0/2^4, \tau_0/2^6$
&$\ep_0/2^6, \tau_0/2^9$\\ \hline
$n^\ep(t=1)$& 2.15E-2& 1.23E-3&6.20E-5&3.88E-6\\
order in time &-& 4.13/3&4.31/3&4.00/3\\
\toprule
\bottomrule
Case II ($\tau=O(\ep)$)
 &$\ep_0=1/2^3, \tau_0=0.1/2^3$
&$\ep_0/2, \tau_0/2$
&$\ep_0/2^2, \tau_0/2^2$
&$\ep_0/2^3, \tau_0/2^3$ \\
\hline
$n^\ep(t=1)$&1.31E-3& 5.12E-4&2.25E-4&1.04E-4\\
order in time &-& 1.35&1.19&1.11\\

\end{tabular*}
{\rule{\temptablewidth}{1pt}}
\end{table}

\section{Conclusion}
A uniformly accurate finite difference method was presented for the
Zakharov system (ZS) with a dimensionless parameter $0<\ep\le1$
which is inversely proportional to the speed of sound.
When $0<\ep\ll1$, i.e. subsonic limit regime,
the solution of ZS propagates highly oscillatory waves in time and/or
rapid outgoing waves in space. Our method was
designed by reformulating ZS into an asymptotic
consistent formulation and adopting an integral approximation
of the oscillating term. Two error bounds were
established by using the energy method and
the limiting equation, respectively,
which depend explicitly on the mesh size $h$ and time step $\tau$ as well as
the parameter $0<\ep\le1$. From the two error bounds, uniform
error estimates were obtained for $0<\ep\le1$.
Numerical results were reported to demonstrate
that the error bounds are sharp.

\section*{Acknowledgements}
This work was partially done while the authors were visiting the Fields
Institute for Research in Mathematical Sciences in Toronto in 2016.


\end{document}